\documentclass[letterpaper,12pt,onecolumn]{IEEEtran}
\usepackage{etoolbox}
\newtoggle{col2}
\togglefalse{col2}
\newtoggle{IEEE}
\togglefalse{IEEE}
\newtoggle{sds}
\toggletrue{sds}

\usepackage{graphicx}
\usepackage{bbm}
\usepackage{dsfont}
\usepackage{amsthm}
\usepackage{amsmath}
\usepackage{amssymb}
\usepackage[mathscr]{eucal}
\usepackage{url}
\usepackage[font=scriptsize,caption=false]{subfig}
\usepackage{color}
\usepackage{cite}

\newcommand{\see}[1]{(see~#1)}

\newcommand{\set}[1]{\left\{ #1 \right\}}
\newcommand{\paren}[1]{\left( #1 \right)}
\newcommand{\abs}[1]{\left| #1 \right|}
\newcommand{\norm}[1]{\left\| #1 \right\|}
\newcommand{\inorm}[1]{\left\| #1 \right\|_i}

\newcommand{\mat}[2]{\left(\begin{array}{#1} #2 \end{array}\right)}
\newcommand{\dt}[1]{\frac{d\,}{d#1}}
\newcommand{\vf}[1]{\frac{\partial\ }{\partial #1}}
\newcommand{\vfof}[2]{\frac{\partial #1}{\partial #2}}
\newcommand{\bd}{\partial}

\newcommand{\rank}[1]{\operatorname{rank}\,#1}
\newcommand{\td}[1]{\widetilde{#1}}
\newcommand{\what}[1]{\widehat{#1}}
\newcommand{\eqcls}[1]{\left[#1\right]}
\newcommand{\qsim}{{\sim}}
\newcommand{\quot}[1]{#1/\qsim}

\newcommand{\eqn}[1]{\eqnal{#1}}
\newcommand{\eqnn}[1]{\eqnaln{#1}}
\newcommand{\eqni}[1]{$#1$}
\newcommand{\eqnal}[1]{\begin{equation*}\begin{aligned}[b] #1 \end{aligned}\end{equation*}}
\newcommand{\eqnaln}[1]{\begin{equation}\begin{aligned} #1 \end{aligned}\end{equation}}
\newcommand{\sm}{\setminus}
\newcommand{\into}{\rightarrow}

\newcommand{\goesto}{\rightarrow}
\newcommand{\inc}{\hookrightarrow}

\newcommand{\ax}[1]{#1^\veps}

\newcommand{\Ad}{\operatorname{Ad}}

\newcommand{\R}{\mathbb{R}}

\newcommand{\N}{\mathbb{N}}

\newcommand{\e}{\mathscr}
\newcommand{\eps}{\epsilon}
\newcommand{\vphi}{\varphi}
\newcommand{\veps}{\varepsilon}
\newcommand{\dx}{\delta x}
\newcommand{\dtheta}{\delta\theta}
\newcommand{\stance}{\operatorname{stance}}
\newcommand{\swing}{\operatorname{swing}}

\newcommand{\rk}{\operatorname{rank}}
\newcommand{\diag}{\operatorname{diag}}
\newcommand{\spec}{\operatorname{spec}}
\newcommand{\specr}[1]{\rho\paren{#1}}

\newcommand{\Pmap}{Poincar\'{e} map}

\iftoggle{IEEE}{
\newcommand{\revi}[1]{{\normalsize{{{#1}}}}}
\newcommand{\revii}[1]{{\normalsize{{{\color{blue}#1}}}}}
}
{
\newcommand{\revi}[1]{{\normalsize{{{#1}}}}}
\newcommand{\revii}[1]{{\normalsize{{{#1}}}}}
}

\newcommand{\sect}[1]{\section{#1}}
\newcommand{\subsect}[1]{\subsection{#1}}

\newtheorem{theorem}{Theorem}
\newtheorem{proposition}{Proposition}
\newtheorem{corollary}{Corollary}
\newtheorem{lemma}{Lemma}
\newtheorem{assumption}{Assumption}

\newtheorem{remark}{Remark}
\newtheorem{example}{Example}
\newtheorem{definition}{Definition}

\newcommand{\defn}[1]{\begin{definition} #1 \end{definition}}

\newcommand{\rem}[1]{\begin{remark} #1 \end{remark}}
\newcommand{\assump}[1]{\begin{assumption} #1 \end{assumption}}
\newcommand{\thm}[1]{\begin{theorem} #1 \end{theorem}}
\newcommand{\lem}[1]{\begin{lemma} #1 \end{lemma}}
\newcommand{\cor}[1]{\begin{corollary} #1 \end{corollary}}
\newcommand{\pf}[1]{\begin{proof} #1 \end{proof}}
\newcommand{\prop}[1]{\begin{proposition} #1 \end{proposition}}

\begin{document}

\title{Model Reduction Near Periodic Orbits \\ of Hybrid Dynamical Systems}

\author{Samuel~A.~Burden\iftoggle{IEEE}{,~\IEEEmembership{Student,~IEEE}}{}, Shai~Revzen\iftoggle{IEEE}{,~\IEEEmembership{Member,~IEEE}}{} and~S.~Shankar~Sastry\iftoggle{IEEE}{,~\IEEEmembership{Fellow,~IEEE}}{}%
\thanks{S.~A.~Burden and S.~S.~Sastry are with the Department of Electrical Engineering and Computer Sciences,
        University of California, Berkeley, CA, USA
        {\tt\small sburden,sastry@eecs.berkeley.edu}}%
\thanks{S.~Revzen is with the Department of Electrical Engineering and Computer Science, University of Michigan,
        Ann Arbor, MI, USA
        {\tt\small shrevzen@eecs.umich.edu}}%
}

\maketitle

\begin{abstract}
We show that, near periodic orbits, a class of hybrid models can be reduced to or approximated by smooth continuous--time dynamical systems.
Specifically, near an exponentially stable periodic orbit undergoing isolated transitions in a hybrid dynamical system, nearby executions generically contract superexponentially to a constant--dimensional subsystem.
Under a non--degeneracy condition on the rank deficiency of the associated Poincar\'{e} map, the contraction occurs in finite time regardless of the stability properties of the orbit.
Hybrid transitions may be removed from the resulting subsystem via a topological quotient that admits a smooth structure to yield an equivalent smooth dynamical system.
We demonstrate reduction of a high--dimensional underactuated mechanical model for terrestrial locomotion, assess structural stability of deadbeat controllers for rhythmic locomotion and manipulation, and derive a normal form for the stability basin of a hybrid oscillator. 
These applications illustrate the utility of our theoretical results for synthesis and analysis of feedback control laws for rhythmic hybrid behavior.
\end{abstract}

\IEEEpeerreviewmaketitle

\section{Introduction}
Rhythmic phenomena are pervasive,
\revi{appearing in physical situations as diverse as legged locomotion~\cite{HolmesEtAl2006}, dexterous manipulation~\cite{BuehlerKoditschek1994}, gene regulation~\cite{GlassPasternack1978}, and electrical power generation~\cite{HiskensReddy2007}}.
The most natural dynamical models for these systems are piecewise--defined or discontinuous owing to intermittent changes in the mechanical contact state of a locomotor or manipulator, or to rapid switches in protein synthesis or constraint activation in a gene or power network.
Such \emph{hybrid} systems generally exhibit dynamical behaviors that are distinct from those of \emph{smooth} systems~\cite{LygerosJohansson2003}.
Restricting our attention to the dynamics near periodic orbits in hybrid dynamical systems, we demonstrate that a class of hybrid models for rhythmic phenomena reduce to classical (smooth) dynamical systems.

Although the results of this paper do not depend on the phenomenology of the physical system under investigation, a principal application domain for this work is terrestrial locomotion.
Numerous architectures have been proposed to explain how animals control their limbs;
for steady--state locomotion, most posit a principle of coordination, synergy, symmetry or synchronization, and there is a surfeit of neurophysiological data to support these hypotheses~\cite{Grillner1985, CohenHolmesRand1982, GolubitskyEtAl1999, TingMacph2005, LiZhang2013}.
Taken together, the empirical evidence suggests that the large number of degrees--of--freedom (DOF) available to a locomotor can collapse during regular motion to a low--dimensional dynamical attractor (a \emph{template}) embedded within a higher--dimensional model (an \emph{anchor}) that respects the locomotor's physiology~\cite{HolmesEtAl2006, FullKoditschek1999}.
We provide a mathematical framework to model this empirically observed dimensionality reduction \revi{in the deterministic setting}.

A stable hybrid periodic orbit provides a natural abstraction for the dynamics of steady--state legged locomotion.  
This widely--adopted approach has generated models of bipedal~\cite{McGeer1990, GrizzleAbba2002, SeyfarthGeyer2003, CollinsRuina2005} and multi--legged~\cite{GhigliazzaAltendorfer2003, SchmittHolmes2000i, KukillayaProctor2009} locomotion as well as control--theoretic techniques for composition~\cite{KlavinsKoditschek2002}, coordination~\cite{HaynesRizzi2012}, and stabilization~\cite{WesterveltGrizzle2003, CarverCowan2009, ShiriaevFreidovich2010}.
In certain cases, it has been possible to embed a low--dimensional abstraction in a higher--dimensional model~\cite{PoulakakisGrizzle2009, AnkaraliSaranli2011}.
Applying these techniques to establish existence of a reduced--order subsystem imposes stringent assumptions on the dynamics of locomotion that are difficult to verify for any particular locomotor.
In contrast, the results of this paper imply that hybrid dynamical systems generically exhibit dimension reduction near periodic orbits solely due to the interaction of the discrete--time switching dynamics with the continuous--time flow.

Under the hypothesis that iterates of the \Pmap\, associated with a periodic orbit in a hybrid dynamical system are eventually constant rank, we construct 
a constant--dimensional invariant subsystem that attracts all nearby trajectories in finite time regardless of the stability properties of the orbit; this appears as Theorem~\ref{thm:exact} of Section~\ref{sec:hds:exact}, below.
Assuming instead that the periodic orbit under investigation is exponentially stable, we show in Theorem~\ref{thm:approx} of Section~\ref{sec:hds:approx} that trajectories \emph{generically} contract superexponentially to a subsystem whose dimension is determined by rank properties of the linearized \Pmap\, at a single point.
The resulting subsystems possess a special structure that we exploit in Theorem~\ref{thm:smooth} to construct a topological quotient that removes the hybrid transitions and admits the structure of a smooth manifold, yielding an equivalent smooth dynamical system.

In Section~\ref{sec:ex} we 
apply these results to reduce the complexity of hybrid models for mechanical systems and analyze rhythmic hybrid control systems.
The example in Section~\ref{sec:hop} demonstrates that reduction can occur spontaneously in mechanical systems undergoing plastic impacts. 
In Section~\ref{sec:poly} we prove that a family of $(3+2n)$--DOF multi--leg models provably reduce to a common 3--DOF mechanical system independent of the number of limbs, $n\in\N$; this demonstrates model reduction in the mechanical component of the class of neuromechanical models considered in~\cite{HolmesEtAl2006, KukillayaProctor2009}.
As further applications, we assess structural stability of deadbeat controllers for rhythmic locomotion and manipulation in Section~\ref{sec:deadbeat} and derive a normal form for the stability basin of a hybrid oscillator in Section~\ref{sec:floq}. 

\sect{Preliminaries}\label{sec:prelim}
We assume familiarity with differential topology and geometry~\cite{Hirsch1976, Lee2002}, and summarize notation and terminology in this section for completeness.

%

If $(X,\norm{\cdot})$ is a Banach space, we let $B_\delta(x)\subset X$ denote the open ball of radius $\delta > 0$ centered at $x\in X$; 
For $X = \R^n$, we may emphasize the dimension $n$ by writing $B^n_\delta(0)\subset\R^n$ for the open $\delta$--ball.
A subset of a topological space is \emph{precompact} if it is open and its closure is compact.
A \emph{neighborhood} of a point $x\in X$ in a topological space $X$ is a connected open subset $U\subset X$ containing $x$.
The \emph{disjoint union} of a collection of sets $\set{S_j}_{j \in J}$ is denoted $\coprod_{j\in J} S_j = \bigcup_{j\in J} S_j \times \set{j}$, a set we endow with the natural piecewise--defined topology.
If $\qsim\subset D\times D$ is an equivalence relation on the topological space $D$, then we let $D/\qsim$ denote the corresponding set of equivalence classes.
There is a natural \emph{quotient projection} $\pi:D\into \quot{D}$ sending $x\in D$ to its equivalence class $\eqcls{x}\in\quot{D}$, and we endow $\quot{D}$ with the (unique) finest topology making $\pi$ continuous~\cite[Appendix~A]{Lee2002}.
Any map $R:G\into D$ defined over a subset $G\subset D$ determines an equivalence relation $\qsim = \set{(x,y)\in D\times D : x\in R^{-1}(y),\ y\in R^{-1}(x),\ \text{or}\  x = y}$.
To be explicit that the equivalence relation is determined by $R$ we denote the quotient space as
\eqn{
\quot{D} = \frac{D}{G \overset{R}{\sim}R(G)}.
}


A \emph{$C^r$ $n$--dimensional manifold} $M$ \emph{with boundary} $\bd M$ is an $n$--dimensional topological manifold covered by an \emph{atlas} of \emph{$C^r$ coordinate charts} $\set{(U_\alpha,\vphi_\alpha)}_{\alpha\in\e{A}}$ where $U_\alpha\subset M$ is open, $\vphi_\alpha:U_\alpha\into H^n$ is a homeomorphism, and $H^n =\set{(y_1,\dots,y_n)\in\R^n : y_n \ge 0}$ is the upper half--space; we write $\dim M = n$.
The charts are $C^r$ in the sense that $\vphi_\alpha\circ\vphi_\beta^{-1}$ is a $C^r$ diffeomorphism over $\vphi_\beta(U_\alpha\cap U_\beta)$ for all pairs $\alpha,\beta\in\e{A}$ for which $U_\alpha\cap U_\beta\ne\emptyset$;
if $r = \infty$ we say $M$ is \emph{smooth}.
The boundary $\bd M\subset M$ contains those points that are mapped to the plane $\set{(y_1,\dots,y_n)\in\R^n : y_n = 0}$ in some chart.
A map $P:M\into N$ is $C^r$ if $M$ and $N$ are $C^r$ manifolds and for every $x\in M$ there is a pair of charts $(U,\vphi),(V,\psi)$ with $x\in U\subset M$ and $P(x)\in V\subset N$ such that the coordinate representation $\td{P} = \psi\circ P\circ\vphi^{-1}$ is a $C^r$ map between subsets of $H^n$.
We let $C^r(M,N)$ denote the normed vector space of $C^r$ maps between $M$ and $N$ endowed with the uniform $C^r$ norm~\cite[Chapter~2]{Hirsch1976}.

Each $x\in M$ has an associated \emph{tangent space} $T_x M$, and the disjoint union of the tangent spaces is the \emph{tangent bundle} $TM = \coprod_{x\in M}T_x M$.
Note that any element in $TM$ may be regarded as a pair $(x,\delta)$ where $x\in M$ and $\delta\in T_x M$, and $TM$ is naturally a smooth $2n$--dimensional manifold.
We let $\e{T}(M)$ denote the set of \emph{smooth vector fields} on $M$, i.e. smooth maps $F:M\into TM$ for which $F(x) = (x,\delta)$ for some $\delta\in T_x M$ and all $x\in M$.
It is a fundamental result that any $F\in\e{T}(M)$ determines an ordinary differential equation in every chart on the manifold that may be solved globally to obtain a \emph{maximal flow} $\phi:\e{F}\into M$ where $\e{F}\subset \R\times M$ is the \emph{maximal flow domain}~\cite[Theorem~17.8]{Lee2002}.  

If $P:M\into N$ is a smooth map between smooth manifolds,
then at each $x\in M$ there is an associated linear map $DP(x):T_xM\into T_{P(x)}N$ called the \emph{pushforward}.
Globally, the pushforward is a smooth map $DP : TM\into TN$; in coordinates, it is the familiar Jacobian matrix.
If $M = X\times Y$ is a product manifold, the pushforward naturally decomposes as $DP = \paren{D_x P,\, D_y P}$ corresponding to derivatives taken with respect to $X$ and $Y$, respectively.
The \emph{rank} of a smooth map $P:M\into N$ at a point $x\in M$ is $\rk DP(x)$.  
If $\rk DP(x) = r$ for all $x\in M$, we simply write $\rk DP \equiv r$. 
If $P$ is furthermore a homeomorphism onto its image, then $P$ is a \emph{smooth embedding}, and the image $P(M)$ is a \emph{smooth embedded submanifold}.
In this case the difference $\dim N - \dim P(M)$ is called the \emph{codimension} of $P(M)$,
and any smooth vector field $F\in\e{T}(M)$ may be pushed forward to a unique smooth vector field $DP(F)\in\e{T}(P(M))$.
A vector field $F\in\e{T}(M)$ is \emph{inward--pointing} at $x\in\bd M$ if for any coordinate chart $(U,\vphi)$ with $x\in U$ the $n$--th coordinate of $D\vphi(F)$ is positive and \emph{outward--pointing} if it is negative.

\sect{Hybrid Dynamical Systems}\label{sec:hds}

We describe a class of hybrid systems useful for modeling physical phenomena in Section~\ref{sec:hdg}, then restrict our attention to the behavior of such systems near periodic orbits in Section~\ref{sec:gamma}.  
It was shown in~\cite{WendelAmes2012} that the Poincar\'{e} map associated with a periodic orbit of a hybrid system is generally not full rank;
we explore the geometric consequences of this rank loss.
Under a non--degeneracy condition on this rank loss we demonstrate in Section~\ref{sec:hds:exact} that the hybrid system possesses an invariant hybrid subsystem to which all nearby trajectories contract in finite time regardless of the stability properties of the orbit.
In Section~\ref{sec:hds:approx} we show that the invariance and contraction of the subsystem hold approximately for any exponentially stable hybrid periodic orbit.
Using tools from differential topology, we remove hybrid transitions from the resulting reduced--order subsystems in Section~\ref{sec:hds:smooth} to yield a continuous--time dynamical system that governs the behavior of the hybrid system near its periodic orbit.

\subsect{Hybrid Differential Geometry}\label{sec:hdg}

For our purposes, it is expedient to define hybrid dynamical systems over a finite disjoint union $M = \coprod_{j\in J} M_j$ where $M_j$ is a connected manifold with boundary for each $j\in J$;
we endow $M$ with the natural (piecewise--defined) topology and smooth structure.
We refer to such spaces as \emph{smooth hybrid manifolds}.
Note that the dimensions of the constituent manifolds are not required to be equal.
Several differential--geometric constructions naturally generalize to such spaces; we prepend the modifier `hybrid' to make it clear when this generalization is invoked.
For instance, the \emph{hybrid tangent bundle} $TM$ is the disjoint union of the tangent bundles $TM_j$,
and the \emph{hybrid boundary} $\bd M$ is the disjoint union of the boundaries $\bd M_j$.

Let $M = \coprod_{j\in J} M_j$ and $N = \coprod_{\ell\in L} N_\ell$ be two hybrid manifolds.
Note that if a map $R:M\into N$ is continuous, then for each $j\in J$ there exists $\ell\in L$ such that $R(M_j)\subset N_\ell$ and hence $R|_{M_j}:M_j\into N_\ell$. 
Using this observation, there is a natural definition of differentiability for continuous maps between hybrid manifolds.
Namely, a map $R:M\into N$ is called \emph{smooth} if $R$ is continuous and $R|_{M_j}:M_j\into N$ is smooth for each $j\in J$.
In this case the \emph{pushforward} $DR:TM\into TN$ is the smooth map defined piecewise as $DR|_{TM_j} = D(R|_{M_j})$ for each $j\in J$.
A smooth map $F:M\into TM$ is called a \emph{vector field} if for all $x\in M$ there exists $v\in T_x M$ such that $F(x) = (x,v)$.

With these preliminaries established, we define the class of hybrid systems considered in this paper.
This is a specialization of \emph{hybrid automata}~\cite{LygerosJohansson2003} that emphasizes the differential--geometric character of hybrid phenomena.

\defn{A \emph{hybrid dynamical system} is specified by a tuple $H = (D,F,G,R)$ where:
\begin{itemize}
\item[$D$] $= \coprod_{j\in J} D_j$ is a smooth hybrid manifold;
\item[$F$] $:D\into TD$ is a smooth vector field;
\item[$G$] $\subset\bd D$ is an open subset of $\bd D$;
\item[$R$] $:G\into D$ is a smooth map.
\end{itemize}
As in~\cite{LygerosJohansson2003}, we call $R$ the \emph{reset map} and $G$ the \emph{guard}.
When we wish to be explicit about the order of smoothness, we will say $H$ is $C^r$ if $D$, $F$, and $R$ are $C^r$ as a manifold, vector field, and map, respectively, for some $r\in\N$.
\label{defn:hds}
}

Roughly speaking, an \emph{execution} of a hybrid dynamical system is determined from an initial condition in $D$ by following the continuous--time dynamics determined by the vector field $F$ until the trajectory reaches the guard $G$, at which point the reset map $R$ is applied to obtain a new initial condition.

%
%

\defn{An \emph{execution} of a hybrid dynamical system $H = (D,F,G,R)$ is a right--continuous function $x:T\into D$ over an interval $T\subset\R$ such that:
\begin{enumerate}
\item if $x$ is continuous at $t\in T$, then $x$ is differentiable at $t$ and $\dt{t} x(t) = F(x(t))$;
\item if $x$ is discontinuous at $t\in T$, then the limit $x(t^-) = \displaystyle\lim_{s\goesto t^-} x(s)$ exists, $x(t^-)\in G$, and $\displaystyle R(x(t^-)) = x(t)$.
\end{enumerate}
\label{defn:exec}
}

If $F$ is tangent to $G$ at $x\in G$, there is a possible ambiguity in determining a trajectory from $x$ since one may either follow the flow of $F$ on $D$ or apply the reset map to obtain a new initial condition $y = R(x)$.

\assump{\label{assump:trans}
$F$ is outward--pointing on $G$.
}

\rem{
The use of time--invariant vector fields and reset maps in Definition~\ref{defn:hds} is without loss of generality in the following sense.
Suppose $D$ is a hybrid manifold, $G\subset\bd D$ is open, and $F:\R\times D\into TD$, $R:\R\times G\into D$ define a time--varying vector field and reset map, respectively.
Define
\eqn{
\revii{\what{D}} = \R\times D,\ \revii{\what{G}} = \R\times G,
}
and let $\revii{\what{F}}:\revii{\what{D}}\into T\revii{\what{D}}$, $\revii{\what{R}}:\revii{\what{G}}\into\revii{\what{D}}$ be defined in the obvious way.
Then $\revii{\what{H}} = (\revii{\what{D}},\revii{\what{F}},\revii{\what{G}},\revii{\what{R}})$ is a hybrid dynamical system in the form of Definition~\ref{defn:hds}.
}

\subsect{Hybrid Periodic Orbits and Hybrid Poincar\'{e} Maps}\label{sec:gamma}

In this paper, we are principally concerned with \emph{periodic} executions of hybrid dynamical systems, which are nonequilibrium trajectories that intersect themselves.

\defn{
An execution $\gamma:T\into D$ is \emph{periodic} if there exists $s\in T$, $\tau > 0$ such that $s + \tau\in T$ and
\eqnn{
\gamma(s) = \gamma(s+\tau).
\label{eqn:gamma}
}
If there is no smaller positive number $\tau$ such that \eqref{eqn:gamma} holds, then $\tau$ is called the \emph{period} of $\gamma$, and we will say $\gamma$ is a \emph{$\tau$--periodic orbit}.
\label{defn:gamma}
}

\rem{
The domain $T$ of a periodic orbit may be taken to be the entire real line, $T = \R$, without loss of generality.
In the sequel we conflate the execution $\gamma:\R\into D$ with its image $\gamma(\R)\subset D$.
}


%
Motivated by the applications in Section~\ref{sec:ex}, we restrict our attention to periodic orbits undergoing \emph{isolated discrete transitions}, i.e. a finite number of discrete transitions that occur at distinct time instants.

\revii{
\assump{\label{assump:isolated}
$\gamma$ undergoes isolated discrete transitions.
}
}

\noindent
In addition to excluding \emph{Zeno} periodic orbits~\cite{OrAmes2011} from our analysis, this assumption enables us to construct Poincar\'{e} maps~\see{\cite{HirschSmale1974, GuckenheimerHolmes1983} for the classical case} associated with 
$\gamma$.
A \emph{Poincar\'{e} map} $P:U\into\Sigma$ is defined over an open subset $U\subset\Sigma$ of an embedded codimension--1 submanifold $\Sigma\subset D$ that intersects the periodic orbit at one point $\set{\xi} = \gamma\cap\Sigma$ by tracing an execution from $x\in U$ forward in time until it intersects $\Sigma$ at $P(x)$.
The submanifold $\Sigma$ is referred to as a \emph{Poincar\'{e} section}.
It is known that this procedure yields a map that is well--defined and smooth near the fixed point $\xi = P(\xi)$~\cite{AizermanGantmacher1958, GrizzleAbba2002, WendelAmes2012, BurdenRevzen2011}. 
Unlike the classical case, Poincar\'{e} maps in hybrid systems need not be full rank.

A straightforward application of Sylvester's inequality~\cite[Appendix~A.5.3]{CallierDesoer1991} shows that the rank of the Poincar\'{e} map is bounded above by the minimum dimension of all hybrid domains. 
More precise bounds are pursued elsewhere~\cite{WendelAmes2012}, but the following Proposition will suffice for the Applications in Section~\ref{sec:ex}.
\prop{If $P:U\into\Sigma$ is a Poincar\'{e} map associated with a periodic orbit $\gamma$, then $\forall x\in U : \rk DP(x) \le \min_{j\in J}\dim D_j - 1$.
\label{prop:rank}}

It is a standard result for continuous--time dynamical systems that the eigenvalues of the linearization of the \Pmap\ at its fixed point---commonly called \emph{Floquet multipliers}---do not depend on the choice of Poincar\'{e} section~\cite[Section~1.5]{GuckenheimerHolmes1983}.
\revi{
This generalizes to the hybrid setting in the sense that there exist similarity transformations relating the non--nilpotent portion of the Jordan forms for linearizations of Poincar\'{e} maps defined over different sections.
Note that, since Proposition~\ref{prop:rank} implies that zero eigenvalues will generally have different algebraic multiplicity for linearized Poincar\'{e} maps obtained from sections located in hybrid domains with different dimensions, we do not expect the nilpotent Jordan blocks for these linear maps to bear any relation to one another.
These observations are summarized in the following Lemma.

\lem{\label{lem:jordan}
If $P:U\into\Sigma$, $\revii{\what{P}}:\revii{\what{U}}\into\revii{\what{\Sigma}}$ are Poincar\'{e} maps associated with a periodic orbit $\gamma$ with fixed points $P(\xi) = \xi$, $\revii{\what{P}}(\revii{\what{\xi}}) = \revii{\what{\xi}}$, 
then $\spec DP(\xi)\sm\set{0} = \spec D\revii{\what{P}}(\revii{\what{\xi}})\sm\set{0}$. 
Moreover, with
\eqn{
J = \mat{cc}{A & 0 \\ 0 & N},\ 
\revii{\what{J}} = \mat{cc}{\revii{\what{A}} & 0 \\ 0 & \revii{\what{N}}}
}
denoting the Jordan canonical forms of $DP(\xi)$ and $D\revii{\what{P}}(\revii{\what{\xi}})$, where 
$0\not\in\spec A\cup\spec\revii{\what{A}}$ 
and
$N$, $\revii{\what{N}}$ are nilpotent, 
we conclude that $A$ is similar to $\revii{\what{A}}$.
}
\pf{
The periodic orbit undergoes a finite number of transitions $k\in\N$, so we may index the corresponding sequence of domains as\footnote{We regard subscripts modulo $k$ so that $D_k \equiv D_0$.} $D_1,\dots,D_k$.
Without loss of generality, assume the $D_j$'s are distinct\footnote{Otherwise we can find $\set{B_j}_{j=1}^k$ such that $B_j\subset D_j$ is open, $\bigcup_{j=1}^k B_j$ contains $\gamma$, and $B_i\cap B_j=\emptyset$ if $i\ne j$, then proceed on $\revii{\what{D}} = \coprod_{j=1}^k B_j$.} and
let $\set{\revii{\xi_j}} = \gamma\cap G\cap\bd D_j$ be the exit point of $\gamma$ in $D_j$.
We wish to construct the Poincar\'{e} map $P_j$ associated with the periodic orbit over a neighborhood of $\xi_j$ in $G$.
For $j \in\set{1,\dots,k}$ let:
\begin{enumerate}
\item[$\phi_j$] $:\e{F}_j\into D_j$ be the maximal flow of $F|_{D_j}$ on $D_j$;
\item[$U_j$] $\subset D_j$ be a neighborhood of $R(\xi_{j-1})$ over which Lemma~\ref{lem:sec} \revii{from Appendix~\ref{sec:tti}} may be applied between $R(\xi_{j-1})\in D_j$ and $G\cap \bd D_j$ to obtain a time--to--impact map $\sigma_j:U_j\into \R$;
\item[$G_j$] $\subset G\cap\bd D_j$ be defined as $G_j = R^{-1}(U_{j+1})$;
\item[$\rho_j$] $:G_j\into G_{j+1}$ be defined by $\rho_j(x) = \phi_{j+1}\paren{\sigma_{j+1}\circ{R(x)},R(x)}$.
\end{enumerate}
The Poincar\'{e} map defined over $G_j$ is obtained formally by iterating the $\rho$'s around the cycle:  
\eqnn{\label{eqn:P} P_j = \rho_{j-1}\circ\cdots\circ \rho_1\circ \rho_k\circ\cdots\circ \rho_j.  }
The neighborhood $\Sigma_j\subset G_j$ of $\xi_j$ over which this map is well--defined is determined by pulling $G_j$ backward around the cycle,
\eqn{\Sigma_j = \paren{\rho_j^{-1}\circ\cdots\circ \rho_k^{-1}\circ \rho_1^{-1}\circ\cdots\circ \rho_{j-1}^{-1}}(G_j),}
and similarly for any iterate of $P_j$.  
Note that $P_j(\xi_j) = \xi_j$ is a fixed point of $P_j$ by construction.
Without loss of generality we assume\footnote{Otherwise we may introduce fictitious guards $\Sigma$ and/or $\revii{\what{\Sigma}}$ near $\gamma$ and repeat the construction.} $\Sigma,\revii{\what{\Sigma}}\subset G$ so that $P = P_j$ and $\revii{\what{P}} = P_i$ for some $i,j\in\set{1,\dots,k}$.
\revii{Refer to Fig.~\ref{fig:pmap} for an illustration.}

\begin{figure}[t]
\centering
\iftoggle{col2}{
\def\svgwidth{1.1\columnwidth} 
\resizebox{1.0\columnwidth}{!}{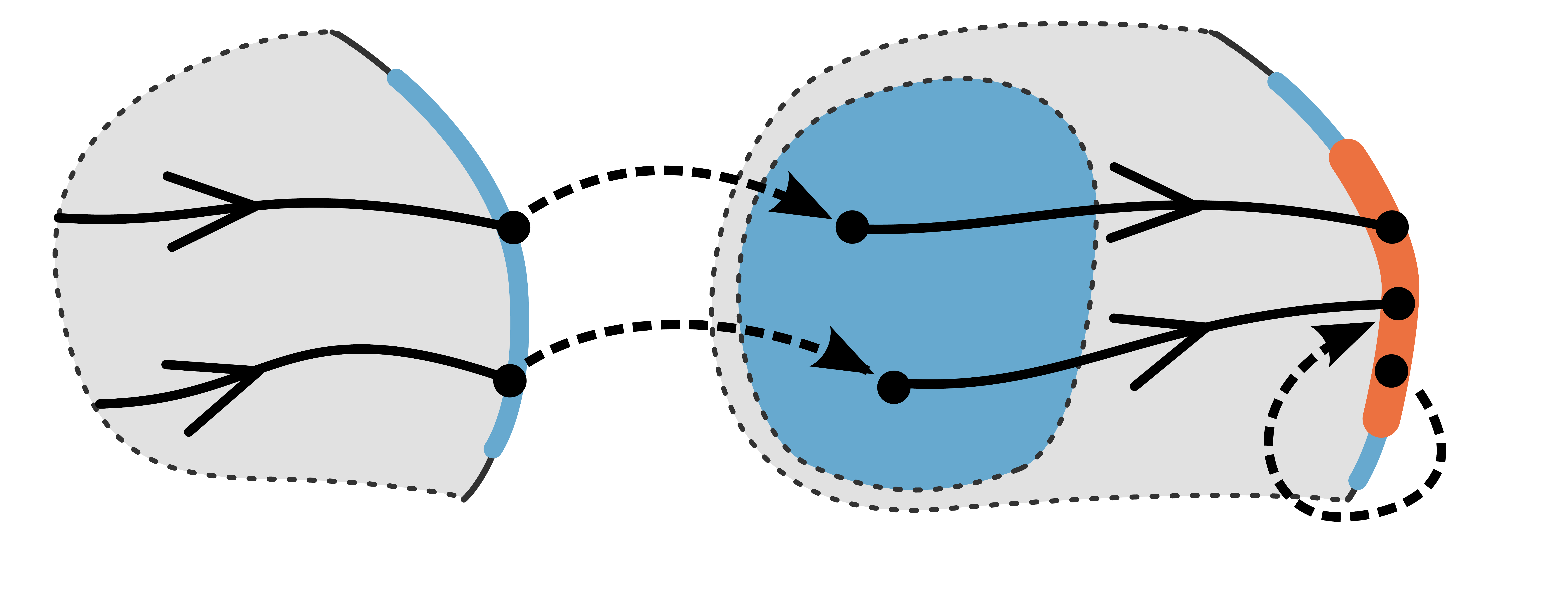}
}
{
\def\svgwidth{0.8\columnwidth} 
\resizebox{0.62\columnwidth}{!}{\input{pmap.pdf_tex}}
}
\caption{%
\revii{\footnotesize
Illustration of constructions used in proofs of Lemma~\ref{lem:jordan}, Theorem~\ref{thm:exact}, and Theorem~\ref{thm:approx}.
For each $j\in J$, the periodic orbit $\gamma$ intersects the guard in domain $D_j$ at $\set{\xi_j} = \gamma\cap G\cap\bd D_j$.
A neighborhood $U_j\subset D_j$ of $R(\xi_{j-1})$  flows via the vector field $F|_{D_j}$ to reach a neighborhood $G_j\subset G\cap\bd D_j$ of $\xi_j$ obtained via $G_j = R^{-1}(U_{j+1})$.
The neighborhood $\Sigma_j\subset G_j$ of $\xi_j$ is chosen sufficiently small to ensure executions initialized in $\Sigma_j$ return to $G_j$ via the {\Pmap} $P_j:\Sigma_j\into G_j$ after one cycle.
}%
\label{fig:pmap}
}
\end{figure}

We proceed by showing that, given a chain of generalized eigenvectors associated with a non--zero eigenvalue of $DP_j(\xi_j)$ for some $j\in\set{1,\dots,k}$, we can construct a chain of generalized eigenvectors associated with $DP_i(\xi_i)$ with the same eigenvalue for each $i\in\set{1,\dots,k}$.
Fix $j\in\set{1,\dots,k}$ and $\lambda\in\spec DP_j(\xi_j)$ with $\lambda\ne 0$.
Suppose $\set{x_j^\ell}_{\ell=1}^m$ is a chain of generalized eigenvectors associated with $\lambda$, i.e. $DP_j(\xi_j)x_j^m = \lambda x_j^m$ and for all $\ell\in\set{1,\dots,m-1}$:
\eqnn{x_j^\ell = (DP_j(\xi_j) - \lambda I) x_j^{\ell+1}.
\label{eqn:eigvec}}
For all $\ell\in\set{1,\dots,m}$, define $x_{j+1}^\ell = D\rho_j(\xi_j) x_j^\ell$ and note $D\rho_j(\xi_j) DP_j(\xi_j) = DP_{j+1}(\xi_{j+1}) D\rho_j(\xi_j)$ by~\eqref{eqn:P}.
Combining this observation with~\eqref{eqn:eigvec} yields 
\eqn{
DP_{j+1}(\xi_{j+1}) x_{j+1}^m 
& = DP_{j+1}(\xi_{j+1}) D\rho_j(\xi_j) x_j^m \\
& = D\rho_j(\xi_j) DP_j(\xi_j) x_j^m \\
& = \lambda D\rho_j(\xi_j) x_j^m = \lambda x_{j+1}^m,
}
so that $\lambda\in\spec DP_{j+1}(\xi_{j+1})$ and for all $\ell\in\set{1,\dots,m-1}$:
\eqn{
x_{j+1}^\ell 
& = D\rho_j(\xi_j) x_j^\ell \\
& = D\rho_j(\xi_j)\paren{DP_j(\xi_j) - \lambda I} x_j^{\ell+1} \\
& = \paren{DP_{j+1}(\xi_{j+1}) D\rho_j(\xi_j) - \lambda D\rho_j(\xi_j)} x_j^{\ell+1} \\
& = (DP_{j+1}(\xi_{j+1}) - \lambda I) x_{j+1}^{\ell+1}.
}
Note that $\set{x_{j+1}^\ell}_{\ell=1}^m$ must be linearly independent since they map to the linearly independent collection $\set{\lambda x_j^\ell}_{\ell=1}^m$ through the composition of linear maps $D\rho_{j-1}(\xi_{j-1})\cdots D\rho_{j+1}(\xi_{j+1})$.
Therefore we conclude $\set{x_{j+1}^\ell}_{\ell=1}^m$ is a chain of generalized eigenvectors for $DP_{j+1}(\xi_{j+1})$ associated with $\lambda$.
Proceeding inductively, for any $i\in\set{1,\dots,k}$ we obtain a corresponding chain for $DP_i(\xi_i)$.
Since the subspace associated with a maximal chain of generalized eigenvectors for a linear map is invariant under the linear map, it follows that the non--nilpotent Jordan blocks of $DP_j(\xi_j)$ must be in one--to--one correspondence with those of $DP_i(\xi_i)$ for any $i\in\set{1,\dots,k}$.
}
}

\subsect{Exact Reduction}\label{sec:hds:exact}

When iterates of the \Pmap\, associated with a periodic orbit of a hybrid dynamical system have constant rank,
executions initialized nearby converge in finite time to a constant--dimensional subsystem.

\thm{[Exact Reduction]\label{thm:exact}
Let $\gamma$ be a periodic orbit that undergoes isolated transitions in a hybrid dynamical system $H = (D,F,G,R)$, 
$P:U\into\Sigma$ a \Pmap\, for $\gamma$, 
\revi{$m = \min_j\dim D_j$},
and suppose there exists a neighborhood $V\subset U$ of $\set{\xi} = \gamma\cap U$ and $r\in\N$ such that $\rk DP^{m}(x) = r$ for all $x\in V$.  
Then there exists an $(r+1)$--dimensional hybrid embedded submanifold $M\subset D$
and a hybrid open set $W\subset D$ for which 
$\gamma\subset M\cap W$ and
trajectories starting in $W$ contract to $M$ in finite time.
}

\pf{
\revi{
We begin in step (i) by applying Lemma~\ref{lem:exact} \revii{from Appendix~\ref{sec:sds:exact}} to construct an $r$--dimensional submanifold $S$ of the Poincar\'{e} section $\Sigma$ that is invariant under the Poincar\'{e} map $P$.
Subsequently, in (ii) we flow $S$ forward in time for one cycle, i.e. until it returns to $\Sigma$, to obtain for each $j\in J$ an $(r+1)$--dimensional submanifold $M_j\subset D_j$ that contains $\gamma\cap D_j$ and is invariant under $F$.
Finally, in (iii) for each $j\in J$ we construct an open set $W_j\subset D_j$ containing $\gamma\cap D_j$ so that the collection $M = \coprod_{j\in J} \revii{M_j}$ attracts all trajectories initialized in the hybrid open set $W = \coprod_{j\in J} \revii{W_j}$ in finite time.
}

(i) 
Applying Lemma~\ref{lem:exact} \revii{from Appendix~\ref{sec:sds:exact}} to $P$, there is a neighborhood $V\subset U$ of $\set{\xi} = \gamma\cap U$ such that $S = P^m(V)$ is an $r$--dimensional embedded submanifold of $U\subset\Sigma$, $P|_S$ maps $S$ diffeomorphically onto $P(S)$, and $P(S)\cap S$ is an open subset of $S$.
Without loss of generality we assume $U\subset G\cap \bd D_1$ and the periodic orbit $\gamma$ passes through each domain once per cycle.
\revi{For notational convenience, for each $j\in J$ we will let $j+1\in J$ denote the subsequent domain visited by $\gamma$ (i.e. we identify $J$ with an additive monoid of integers modulo $\abs{J}$).}
Set $\revii{\set{\xi_1}} = \gamma\cap G\cap\bd D_1$, let $U_2\subset D_2$ be a neighborhood of $R(\revii{\xi_1})$ over which Lemma~\ref{lem:sec} \revii{from Appendix~\ref{sec:tti}} may be applied to construct a time--to--impact map $\sigma_2:U_2\into\R$, let $G_1 = R^{-1}(U_2)$ be a neighborhood of $\xi_1$ in $G\cap \bd D_1$, and let $\phi_1:\e{F}_1\into D_1$ the maximal flow of $F|_{D_1}$ on $D_1$.
Proceed inductively forward around the cycle 
to construct, for each $j\in J$: 
the exit point $\revii{\set{\xi_j}} = \gamma\cap G\cap \bd D_j$; 
time--to--impact map $\sigma_j:U_j\into\R$ over a neighborhood $U_j\subset D_j$ containing $R(\revii{\xi_{j-1}})$; 
a neighborhood $G_j = R^{-1}(U_{j+1})\subset G\cap \bd D_j$ containing $\revii{\xi_j}$;
and the maximal flow $\phi_j:\e{F}_j\into D_j$ of $F|_{D_j}$ on $D_j$.
\revii{Refer to Fig.~\ref{fig:pmap} for an illustration of this construction.}

(ii)
By flowing $S$ forward through one cycle, for each $j\in J$ we will construct a submanifold $M_j\subset D_j$ that is diffeomorphic to $[0,1]\times\R^r$.
Observe that, since $P|_S$ is a diffeomorphism, with $S_1 = S\cap G_1$ we have that the restriction $R|_{S_1}$ is a diffeomorphism onto its image and $F|_{R(S_1)}$ is nowhere tangent to $R(S_1)$.
Let $M_2\subset D_2$ be the embedded submanifold obtained by flowing $R(S_1)$ to $G\cap \bd D_2$, and let $S_2 = M_2\cap G_2$; observe that $S_2$ is diffeomorphic to $S_1$, $M_2$ is diffeomorphic to $[0,1]\times S_2$, and $F|_{D_2}$ is tangent to $M_2$.
Proceed inductively forward around the cycle to construct, for each $j\in J$, an embedded submanifold $S_j\subset G_j$ diffeomorphic to $S_1$ and a submanifold $M_j\subset D_j$ diffeomorphic to $[0,1]\times S_j$ such that $F|_{D_j}$ is tangent to $M_j$.
Note that $S_1$ is diffeomorphic to the $r$--dimensional manifold $\R^r$, so $\dim M_j = r+1$ for each $j\in J$.
The subsystem $M = \coprod_{j\in J} M_j\subset D$ contains $\gamma$, is invariant under the continuous flow by construction, and is invariant under the reset map in the sense that $R^{-1}(M)\cap M\subset G\cap M$ is open.

(iii)
Finally, let $W_1 = \phi_j^{-1}(\R\times V)\subset D_1$ be the open set that flows into $V$, where $S = P^m(V)$ was defined in step (i).
Let $W_{\abs{J}} = \phi_{\abs{J}}^{-1}(R^{-1}(W_1))\subset D_{\abs{J}}$ be the open set that flows into $W_1$ where $\abs{J}$ denotes the number of elements in $J$. 
Proceed inductively backward around the cycle to construct, for each $j\in J$, an open set $W_j\subset D_j$ that flows into $S$ in finite time.
Then the hybrid open set $W = \coprod_{j\in J} W_j\subset D$ contains $\gamma$ and all executions initialized in $W$ flow into $S\subset M$ in finite time.
}

Since $M$ is invariant under the continuous dynamics ($F|_M$ is tangent to $M$) and the discrete dynamics ($R(G\cap M)\subset M$), it determines a subsystem that governs the stability of $\gamma$ in $H$.

\cor{\label{cor:exact:subsys}
$H|_M = (M,F|_M,G\cap M,R|_{G\cap M})$ is a hybrid dynamical system with periodic orbit $\gamma$.
}

\cor{\label{cor:stab}
The periodic orbit $\gamma$ is Lyapunov (resp. asymptotically, exponentially) stable in $H$ if and only if $\gamma$ is Lyapunov (resp. asymptotically, exponentially) stable in $H|_M$.
}

When the rank at the fixed point $\xi = P(\xi)$ achieves the upper bound stipulated by Proposition~\ref{prop:rank}, the following Corollary ensures that $DP^m$ is constant rank (and hence Theorem~\ref{thm:exact} may be applied).
This is important since it is possible to compute a lower bound for $\rk DP^m(\xi)$ via numerical simulation~\cite{BurdenGonzalezVasudevan2013}.

\cor{If $\rk DP^{m}(\xi) = \min_{j\in J}\dim D_j - 1 = m - 1$, then there exists an open set $V\subset U$ containing $\xi$ such that $\rk DP^{m}(x) = m - 1$ for all $x\in V$.
Thus the hypotheses of Theorem~\ref{thm:exact} are satisfied with $r = m - 1$.
\label{cor:exact}}

\revi{
If the Poincar\'{e} map attains the same constant rank $r$ for two subsequent iterates, it is not necessary to continue up to iterate $m = \min_j\dim D_j$ before checking the hypotheses of Theorem~\ref{thm:exact}.

\cor{\label{cor:exact:iterate}
If there exists a neighborhood $W\subset U$ of $\xi$ and $k,r\in\N$ such that $\rank DP^k(x) = r$ for all $x\in W$ and $\rank DP^{k+1}(\xi) = \rank DP^k(\xi)$, then there exists a neighborhood $V\subset W$ of $\xi$ such that $\rank DP^{m}(x) = r$ for all $x\in V$.
Thus the hypotheses of Theorem~\ref{thm:exact} are satisfied with $r = \rank DP^k(\xi)$.
}

}

The choice of Poincar\'{e} section in Theorem~\ref{thm:exact} is irrelevant in the sense that the \Pmap\, $\td{P}:\td{U}\into\td{\Sigma}$ defined over any other Poincar\'{e} section $\td{\Sigma}$ will be constant--rank in a neighborhood $\td{V}\subset\td{U}$ of its fixed point $\set{\td{\xi}} = \gamma\cap\td{\Sigma}$, as the following Corollary shows; this follows directly from Lemma 4 in~\cite{BurdenRevzen2011}.

\cor{\label{cor:pmap}
Under the hypotheses of Theorem~\ref{thm:exact},
if $\td{P}:\td{U}\into\td{\Sigma}$ is any other \Pmap\, for $\gamma$ with fixed point $\td{\xi} = \td{P}(\td{\xi})$, then there exists an open subset $\td{V}\subset\td{U}$ containing $\td{\xi}$ such that $\rk D\td{P}^m(x) = r$ for all $x\in\td{V}$.
Thus the hypotheses of Theorem~\ref{thm:exact} are satisfied for $\td{P}$ with $r = \rank D\td{P}^m(\xi)$.
}

\subsect{Approximate Reduction}\label{sec:hds:approx}

Restricting our attention to exponentially stable periodic orbits, 
we find that a hybrid system generically contracts superexponentially to a constant--dimensional subsystem near a periodic orbit.

\thm{[Approximate Reduction]\label{thm:approx}
Let $\gamma$ be an exponentially stable periodic orbit undergoing isolated transitions in a hybrid dynamical system $H = (D,F,G,R)$, 
$P:U\into\Sigma$ a \Pmap\, for $\gamma$ with fixed point $\set{\xi} = \gamma\cap\Sigma$, 
\revi{$m = \min_j\dim D_j$}, 
and $r = \rk DP^{m}(\xi)$.
Then there exists an $(r+1)$--dimensional hybrid embedded submanifold $M\subset D$
such that for any $\veps > 0$ there exists a hybrid open set $\ax{W}\subset D$ for which 
$\gamma\subset M\cap\ax{W}$ and 
the distance from trajectories starting in $\ax{W}$ to $M$ contracts by $\veps$ each cycle.
}

\pf{
\revi{
We begin with an overview of the proof strategy.
First (i), for each $j\in J$ we construct a Poincar\'{e} map $P_j$ over a Poincar\'{e} section $\Sigma_j\subset G\cap\bd D_j$ and apply Lemma~\ref{lem:approx} \revii{from Appendix~\ref{sec:sds:approx}} to obtain a change--of--coordinates in which $P_j$ splits into two components: a linear map that only depends on the first $r$ coordinates and a nonlinear map whose linearization is nilpotent at the fixed point of $P_j$.
Second (ii), for each $j\in J$ we construct an $r$--dimensional submanifold $S_j\subset\Sigma_j$ such that $R|_{S_j}$ is a diffeomorphism near the fixed point of $P_j$. 
We subsequently flow the image $R(S_j)$ forward until it impacts the guard to construct an $(r+1)$--dimensional submanifold $M_{j+1}\subset D_{j+1}$ that contains $\gamma\cap D_{j+1}$ and is invariant under $F$.
Third (iii), for each $j\in J$ we apply Lemma~\ref{lem:superstable} \revii{from Appendix~\ref{sec:superstable}} to construct a distance metric on an open set $W_j\subset D_j$ containing $\gamma\cap D_j$ with respect to which executions contract superexponentially toward $M_j$.
}

(i)
Without loss of generality we assume $U\subset G\cap \bd D_1$ the periodic orbit $\gamma$ passes through each domain once per cycle.
\revi{As in the proof of Theorem~\ref{thm:exact}, for each $j\in J$ we will let $j+1\in J$ denote the subsequent domain visited by $\gamma$ (i.e. we identify $J$ with an additive monoid of integers modulo $\abs{J}$).}
For each $j\in J$ let $P_j:U_j\into\Sigma_j$ be a \Pmap\, for $\gamma$ defined over $U_j\subset\Sigma_j\subset G\cap\bd D_j$, and let $\revii{\set{\xi_j}} = \gamma\cap G\cap\bd D_j$ be the exit point of $\gamma$ in $D_j$.
\revii{Refer to Fig.~\ref{fig:pmap} for an illustration of this construction.}
\revii{Lemma~\ref{lem:jordan} implies that $\rk DP_j^m(\xi_j) = r$ for all $j\in J$.}
Applying Lemma~\ref{lem:approx} \revii{from Appendix~\ref{sec:sds:approx}} implies that for each $j\in J$ there exists an open set $V_j\subset U_j$ containing $\xi_j$ and a $C^1$ diffeomorphism $\vphi_j:V_j\into\R^{n_j-1}$ where $n_j = \dim D_j$ such that $\vphi_j(\xi_j) = 0$ and the coordinate representation 
$\td{P}_j = \vphi_j\circ P_j\circ\vphi_j^{-1}$ of $P_j$ has the form 
$\td{P}_j(z_j,\zeta_j) = \paren{A_j z_j,\ S_j(z_j,\zeta_j)}$
where $z_j\in\R^r$, $\zeta_j\in\R^{n_j-1-r}$, 
$A_j\in\R^{r\times r}$ is invertible, 
$S_j(0,0) = 0$,
and $D_{\zeta_j} S_j(0,0)$ is nilpotent.
For each $j\in J$, let $\Pi_j:V_j\into G$ be a smooth map defined as follows.
Given $x\in V_j$, write $(z_x,\zeta_x) = \vphi_j(x)\in\R^r\times\R^{n_j-r-1}$ and let $\Pi_j(x) = \vphi_j^{-1}(z_x,0)$.

(ii)
Fix $j\in J$ and let $N_j = \vphi_j^{-1}\paren{\R^r\times\set{0}}\subset V_j$, an $r$--dimensional embedded submanifold tangent to the non--nilpotent eigendirections of $DP_j^m(\xi_j)$.
Observe that $DR|_{G\cap N_j}(\xi_j)$ has rank $r = \dim N_j$, hence by the Inverse Function Theorem~\cite[Theorem~7.10]{Lee2002} there is a neighborhood $S_j\subset N_j$ containing $\xi_j$ such that $R|_{S_j}:S_j\into D$ is a diffeomorphism onto its image $R(S_j)\subset D_{j+1}$.
Furthermore, since $\rk DP_j^m(\xi_j) = r$, the vector field is transverse to $R(S_j)$ at $\xi_j$, i.e. $F(R(\xi_j))\not\in T_{R(\xi_j)} R(S_j)$, and we assume $S_j$ was chosen small enough so that $F$ is transverse along all of $R(S_j)$.
Let $M_{j+1}\subset D_{j+1}$ be the embedded submanifold obtained by flowing $R(S_j)$ forward to $G$; note that $M_{j+1}$ is diffeomorphic to $[0,1]\times\R^r$.
Observe that $M = \coprod_{j\in J}M_j$ is invariant under the continuous flow (i.e. $F|_M$ is tangent to $M$) and approximately invariant under the reset map in the sense that $DR|_{G\cap M}$ is tangent to $M$ on $\gamma$: for all $j\in J$ and $\delta\in T_{\xi_j} (G\cap M)$ we have $DR|_{G\cap M}(\xi_j)\delta\in T_{R(\xi_j)}M$.
Observe that $R\circ\Pi_j|_{G\cap M_j}:G\cap M_j\into M_{j+1}$ is a diffeomorphism onto its image.


(iii)
Fix $\veps > 0$ and apply the construction in the proof of Lemma~\ref{lem:superstable} \revii{from Appendix~\ref{sec:superstable}} to obtain a radius $\delta > 0$ and for each $j\in J$ a norm $\ax{\norm{\cdot}}_j:\R^{n_j-1}\into\R$ such that the nonlinearity $\td{P}_j(z_j,\zeta_j) - \paren{A_j z_j,0}$ contracts exponentially fast with rate $\veps$ on $B^{n_j-1}_\delta(0)\subset\R^{n_j-1}$ as measured by $\ax{\norm{\cdot}}_j$.
For each $j\in J$ 
define $\ax{V}_j = \vphi_j^{-1}(B^{n_j-1}_\delta(0))\subset G\cap \bd D_j$, 
let $\phi_j:\e{F}_j\into D_j$ denote the maximal flow of $F|_{D_j}$ on $D_j$, 
and let $\ax{W}_j = \phi_j^{-1}(\R\times\ax{V}_j)\subset D_j$ be the (open) set of points that flow into $\ax{V}_j$. 
Since $\phi_j$ is the flow of a smooth vector field transverse to $\ax{V}_j$, any $x\in\ax{W}_j$ can be written uniquely as $x = \phi_j(t_x,v_x)$ for some $t_x \le 0$ and $v_x\in\ax{V}_j$.
Using this representation, we endow $\ax{W}_j$ with a distance metric $\ax{d}_j:\ax{W}_j\times\ax{W}_j\into\R$ by defining $\ax{d}_j(x,y) = \abs{t_x - t_y} + \ax{\norm{\vphi_j(v_x) - \vphi_j(v_y)}_j}$.
Observe that the exponential contraction of $\td{P}_j$ at rate $\veps$ in $\ax{\norm{\cdot}}_j$ to $\vphi_j(M_j\cap G)$ implies exponential contraction of executions initialized in $\ax{W}_j$ at rate $\veps$ to $M$ in $\ax{d}_j$.
%

Finally, let $\ax{W} = \coprod_{j\in J}\ax{W}_j$ and $\ax{M} = M\cap\ax{W}$.
Define a smooth hybrid map $\ax{\Pi}:G\cap\ax{W}\into G$ piecewise
for each $j\in J$ by observing that $G\cap\ax{W}_j\subset V_j$ and 
letting $\ax{\Pi}(x) = \Pi_j(x)$ for all $x\in G\cap\ax{W}_j$.
}

\cor{\label{cor:approx:subsys}
Letting $\ax{M} = M\cap\ax{W}$,
the collection $H|_{\ax{M}} = (\ax{M},F|_{\ax{M}},G\cap \ax{M},R\circ\ax{\Pi}|_{G\cap\ax{M}})$ is a $C^1$ hybrid dynamical system with periodic orbit $\gamma$, where $\ax{\Pi}:G\cap\ax{W}\into G$ is the smooth hybrid map constructed in the proof of Theorem~\ref{thm:approx}.
}

Although the submanifold $M\subset D$ is invariant under the continuous dynamics of $H$ in the sense that $F|_M$ is tangent to $M$, the reset map must be modified to ensure $M$ is invariant under the discrete dynamics.
However, since $DR|_{G\cap\ax{M}} = D\paren{R\cap\ax{\Pi}}|_{G\cap\ax{M}}$, the map $\Pi$ does not affect $R$ to first order.

\rem{
We emphasize that hypothesis on the rank of the \Pmap\, $P:U\into\Sigma$ in Theorem~\ref{thm:approx} 
($\rk DP^m(\xi) = r$ at the point $\set{\xi} = \gamma\cap\Sigma$) 
is weaker than the hypothesis in Theorem~\ref{thm:exact} 
($\rk DP^m(x) = r$ for all $x$ in an open set $V\subset U$).
In particular, 
approximating the rank over an uncountably infinite set typically involves estimates on higher--order derivatives of $P^m$.
}

\revi{
If the rank is constant for two subsequent iterates of the linearized Poincar\'{e} map, then the rank is constant for all subsequent iterates, including iterate $m = \min_j\dim D_j$.

\cor{\label{cor:approx:iterate}
If there exist $k\in\N$ such that $\rank DP^k(\xi) = \rank DP^{k+1}(\xi)$, then $\rank DP^m(\xi) = \rank DP^k(\xi)$.
Thus the hypotheses of Theorem~\ref{thm:approx} are satisfied with $r = \rank DP^k(\xi)$.
}
}

\subsect{Smoothing}\label{sec:hds:smooth}

The subsystems yielded by Theorems~\ref{thm:exact}~and~\ref{thm:approx} on \emph{exact} and \emph{approximate} reduction share important properties:
the constituent manifolds have the same dimension; 
the reset map is a hybrid diffeomorphism between disjoint portions of the boundary; 
and the vector field points inward along the range of the reset map. 
Under these conditions, we can globally \emph{smooth} the hybrid transitions using techniques from differential topology to obtain a single continuous--time dynamical system.
Executions of the hybrid (sub)system are preserved as integral curves of the continuous--time system.
This provides a smooth $n$--dimensional generalization of the \emph{hybrifold} construction in~\cite{SimicJohansson2005}, 
the \emph{phase space} constructed in~\cite{Schatzman1998} to analyze mechanical impact,
\revi{as well as the change--of--coordinates constructed in~\cite[\S3.1.1]{De2010} to simplify analysis of juggling}.

\thm{[Smoothing]\label{thm:smooth}
Let $H = (M,F,G,R)$ be a hybrid dynamical system with $M = \coprod_{j\in J}M_j$.  
Suppose 
$\dim M_j = n$ for all $j\in J$, 
$R(G)\subset\bd M$, 
$\bd M = G\coprod R(G)$,  
$R$ is a hybrid diffeomorphism onto its image, 
and $F$ is inward--pointing along $R(G)$. 
Then the topological quotient
$\td{M} = \frac{M}{G \overset{R}{\sim}R(G)}$
may be endowed with the structure of a smooth manifold such that:
\begin{enumerate}
\item the quotient projection $\pi:M\into \td{M}$ restricts to a smooth embedding $\pi|_{M_j}:M_j\into\td{M}$ for each $j\in J$;
\item there is a smooth vector field $\td{F}\in\e{T}(\td{M})$ such that any execution $x:T\into M$ of $H$ descends to an integral curve of $\td{F}$ on $\td{M}$ via $\pi:M\into\td{M}$:
\eqnal{
\forall t\in T : \dt{t}\pi\circ x(t) = \td{F}\paren{\pi\circ x(t)}.
}
\end{enumerate}
}

\pf{
Let $S\subset G\cap M_i$ be a connected component in some domain $i\in J$, and let $k\in J$ be the index for which $R(S)\subset M_k$.
The hypotheses of this Theorem together with Assumption~\ref{assump:trans} ensure
Lemma~\ref{lem:smooth} \revii{from Appendix~\ref{sec:sds:smooth}} may be applied to attach $M_i$ to $M_k$ to yield a new smooth manifold $\td{M}_{ik}$.
The hybrid system defined over the 
\revi{domain $\coprod\set{\td{M}_{ik}}\cup\set{M_j : j\in J\sm\set{i,k}}$} 
and guard $G\sm S$ satisfies the hypotheses of this Theorem, hence we may inductively attach domains on each connected component that remains in $G\sm S$.
This yields a smooth manifold $\td{M}$ and vector field $\td{F}\in\e{T}(\td{M})$ with the required properties.
}

\begin{figure*}[t]
\centering
\subfloat[$H = (D,F,G,R)$]{
\iftoggle{col2}{
\def\svgwidth{0.8\columnwidth} 
\resizebox{0.62\columnwidth}{!}{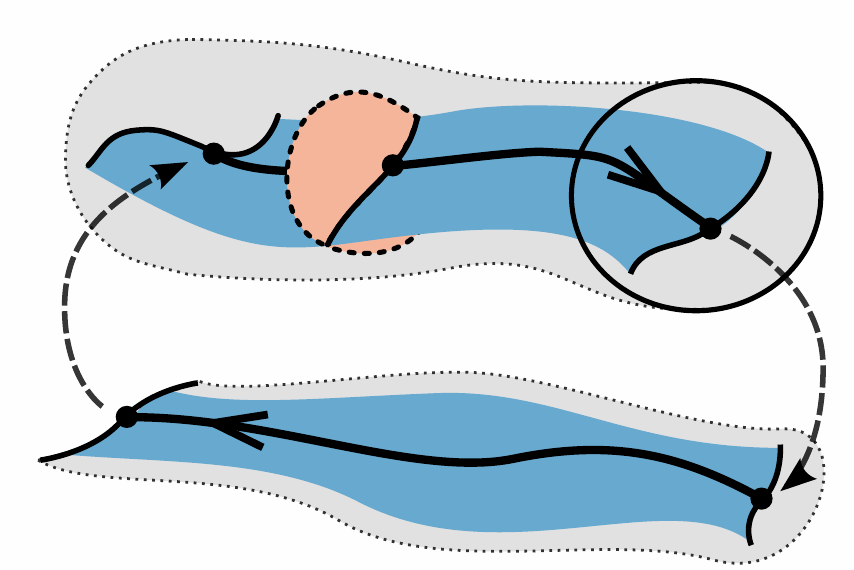}
}{
\def\svgwidth{0.4\columnwidth} 
\resizebox{0.31\columnwidth}{!}{\input{res.pdf_tex}}
}
}
\subfloat[$H|_M = (M,F|_M,G\cap M,R|_{G\cap M})$]{
\iftoggle{col2}{
\def\svgwidth{0.8\columnwidth}
\resizebox{0.62\columnwidth}{!}{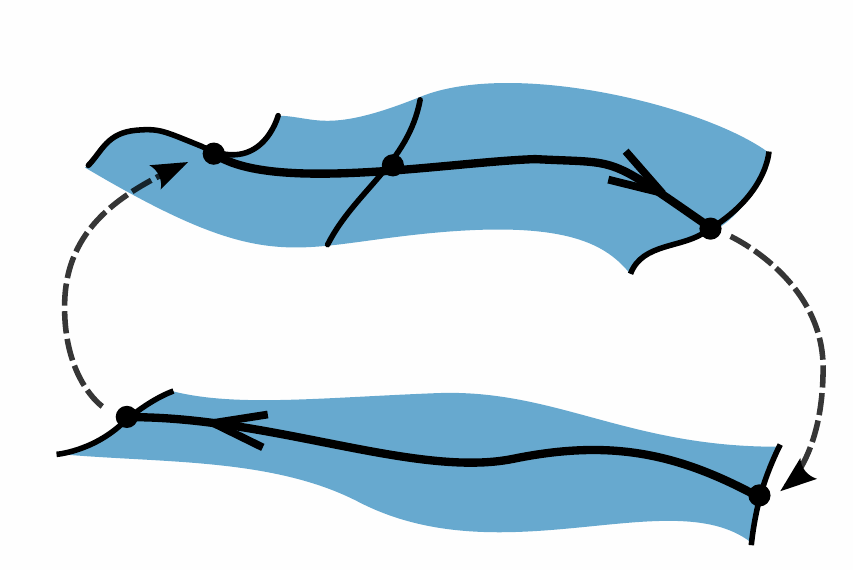}
}{
\def\svgwidth{0.4\columnwidth}
\resizebox{0.31\columnwidth}{!}{\input{reduced.pdf_tex}}
}
}
\subfloat[$(\td{M},\td{F})$]{
\iftoggle{col2}{
\def\svgwidth{0.8\columnwidth}
\resizebox{0.62\columnwidth}{!}{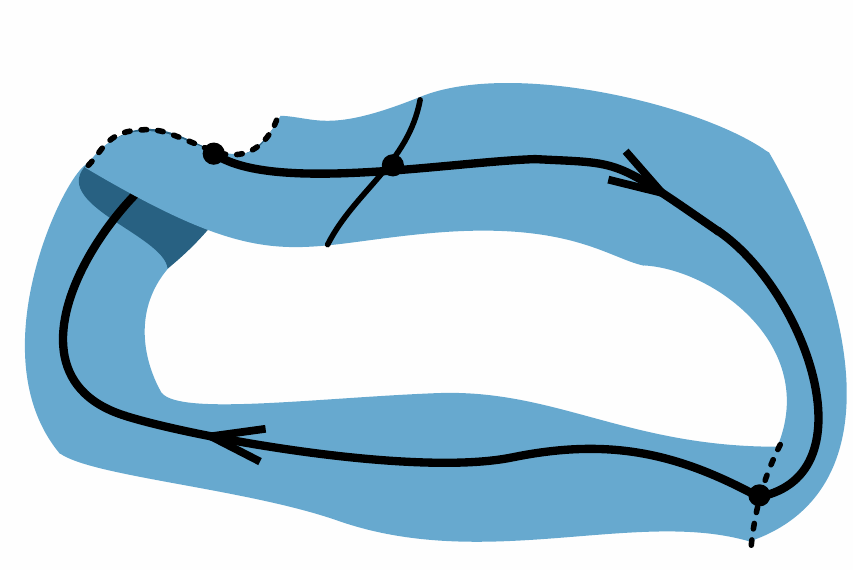}
}{
\def\svgwidth{0.4\columnwidth}
\resizebox{0.31\columnwidth}{!}{\input{smoothed.pdf_tex}}
}
}
\caption{
(a) Applying Theorem~\ref{thm:exact} (Exact Reduction) to a hybrid dynamical system $H = (D,F,G,R)$ containing a periodic orbit $\gamma$ with associated Poincar\'{e} map $P:U\into\Sigma$
yields an invariant subsystem $M = \coprod_{j\in J}M_j$; nearby trajectories contract to $M$ in finite time.
(b) The subsystem may be extracted to yield a hybrid dynamical system $H|_M$.
(c) The hybrid system $H|_M$ may subsequently be smoothed via Theorem~\ref{thm:smooth} (Smoothing) to yield a continuous-time dynamical system $(\td{M},\td{F})$.
Application of Theorem~\ref{thm:smooth} to the subsystem from Theorem~\ref{thm:approx} (Approximate Reduction) is illustrated by replacing $H|_M$ 
with~$H|_{\ax{M}}$.
}
\label{fig:smooth}
\end{figure*}

\rem{
As illustrated in Fig.~\ref{fig:smooth}, Theorem~\ref{thm:smooth} is applicable to the subsystems $H|_M$, $H|_{\ax{M}}$ that emerge as a consequence of the Corollaries to Theorems~\ref{thm:exact}~and~\ref{thm:approx}, respectively.
Thus a class of hybrid models for periodic phenomena may be reduced (exactly or approximately) to smooth dynamical systems.
}

\sect{Applications}\label{sec:ex}

The Theorems of Section~\ref{sec:hds} apply directly to autonomous hybrid dynamical systems; in Section~\ref{sec:hop} we demonstrate that reduction to a smooth subsystem can occur spontaneously in a mechanical system undergoing intermittent impacts.
The results are also applicable to systems with control inputs; in Section~\ref{sec:poly} we synthesize a state--feedback control law that reduces a family of multi--leg models for lateral--plane locomotion to a common low--dimensional subsystem,
and in Section~\ref{sec:deadbeat} we analyze the structural stability of event--triggered deadbeat control laws for locomotion.
Finally, the reduction of hybrid dynamics to a smooth subsystem provides a route through which tools from classical dynamical systems theory can be generalized to the hybrid setting; in Section~\ref{sec:floq} we extend a normal form for limit cycles.

\subsect{Spontaneous Reduction in a Vertical Hopper}\label{sec:hop}

\begin{figure}
\begin{center}
\includegraphics[height=5.0cm]{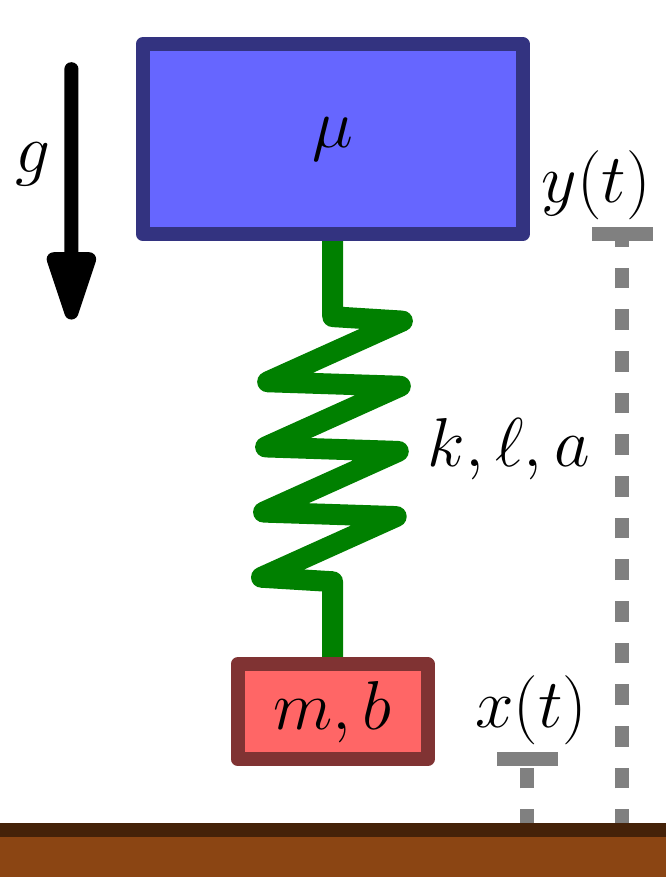}
\caption{
Schematic of vertical hopper.
Two masses $m$ and $\mu$, constrained to move vertically above a ground plane in a gravitational field with magnitude $g$, are connected by a linear spring with stiffness $k$ and nominal length $\ell$.
The lower mass experiences viscous drag proportional to velocity with constant $b$ when it is in the air, and impacts plastically with the ground (i.e. it is not permitted to penetrate the ground and its velocity is instantaneously set to zero whenever a collision occurs).
When the lower mass is in contact with the ground, the spring stiffness is multiplied by a constant $a > 1$.
}
\label{fig:hop}
\end{center}
\end{figure}

\revi{
In this section, we apply Theorem~\ref{thm:exact} (Exact Reduction) to the \emph{vertical hopper} example shown in Fig.~\ref{fig:hop}. 
This system evolves through an \emph{aerial} mode and a \emph{ground} mode.
In the aerial mode, the lower mass moves freely at or above the ground height.
Transition to the ground mode occurs when the lower mass reaches the ground height with negative velocity, where it undergoes a perfectly plastic impact (i.e. its velocity is instantaneously set to zero).
In the ground mode, the lower mass remains stationary.
Transition to the aerial mode occurs when the aerial mode force allows the mass to lift off.
We now formulate this model in the hybrid dynamical system framework of Definition~\ref{defn:hds}. 

The aerial mode $D_a$ (see Fig.~\ref{fig:hop} for notation) consists of
\eqni{
(y,\dot{y},x,\dot{x})\in D_a = T\R\times T\R_{\ge 0},
}
and the vector field $F|_{D_a}$ is given by 
$\mu\ddot{y} = k(\ell - (y-x)) - \mu g$,
$m\ddot{x} = - k(\ell - (y-x)) - b\dot{x} - m g$.
The boundary 
\eqni{
\bd D_a = \set{(y,\dot{y},x,\dot{x})\in D_a : x = 0}
}
contains the states where the lower mass has just impacted the ground,
and a hybrid transition occurs on the subset 
\eqni{
G_a = \set{(y,\dot{y},0,\dot{x})\in\bd D_a : \dot{x} < 0}
}
of the boundary $D_a$ where the lower mass has negative velocity.
The state is reinitialized in the ground mode via $R|_{G_a} : G_a \into D_g$ defined by $R|_{G_a}(y,\dot{y},0,\dot{x}) = (y,\dot{y})$.
In the ground mode $D_g = \set{(y,\dot{y})\in T\R : - k(\ell - y) \le m g}$, the boundary consists of the set of configurations where the force in the aerial mode allows the lower mass to lift off,
$\bd D_g = \set{(y,\dot{y})\in D_g : - k(\ell - y) = m g}$,
and the vector field $F|_{D_g}$ is given by $\mu\ddot{y} = ak(\ell - y) - \mu g$.
A hybrid transition occurs when the forces balance and will instantaneously increase to pull the mass off the ground,
$G_g = \set{(y,\dot{y})\in\bd D_g : \dot{y}(t) > 0}$,
and the state is reset via $R|_{G_g} : G_g \into D_a$ defined by $R|_{G_g}(y,\dot{y}) = (y,\dot{y},0,0)$.
This defines a hybrid dynamical system $(D,F,G,R)$ where
\eqn{
D = D_a\coprod D_g,\
F\in\e{T}(D),\
G = G_a\coprod G_g,\ 
R:G\into D. 
}
}

With parameters $(m,\mu,k,b,\ell,a,g) = (1,3,10,5,2,2,2)$, \revi{numerical simulations suggest} the vertical hopper possesses a stable periodic orbit $\gamma = (y^*,\dot{y}^*,x^*,\dot{x}^*)$ to which nearby trajectories $(y,\dot{y},x,\dot{x})$ converge asymptotically.
Choosing a Poincar\'{e} section $\Sigma$ in the ground domain $D_g$ at mid-stance, $\Sigma = \set{(y,\dot{y}) : \dot{y} = 0}\subset D_g$, we find numerically\footnote{
For numerical simulations, we use a recently--developed algorithm~\cite{BurdenGonzalezVasudevan2013} with step size $h = 1\times 10^{-2}$ and relaxation parameter $\eps = 1\times 10^{-10}$.  
}
that the hopper possesses a stable periodic orbit $\gamma$ that intersects the Poincar\'{e} section at $\gamma\cap\Sigma = \set{\xi}$ where $\xi = (y,\dot{y}) \approx (0.94,0.00)$.
\revi{Using finite differences}, we determine that the linearization $DP$ of the associated scalar--valued Poincar\'{e} map $P:\Sigma\into\Sigma$ has eigenvalue $\spec DP(\xi)\approx 0.57$ at the fixed point $P(\xi) = \xi$.
The rank of the Poincar\'{e} map $P$ attains the upper bound of Proposition~\ref{prop:rank}, hence Corollary~\ref{cor:exact} implies the rank hypothesis of Theorem~\ref{thm:exact} (Exact Reduction) is satisfied.
Thus the dynamics of the hopper collapse to a one degree--of--freedom mechanical system after a single hop.
\revi{Geometrically, the portion of the reduced subsystem in each domain is diffeomorphic to $[0,1]\times\R$.}
Algebraically, the constraint that activates when the lower mass impacts the ground transfers to the aerial mode where no such physical constraint exists: \revi{the lower mass state $(x,\dot{x})$ is uniquely determined by the upper mass state $(y,\dot{y})$ for all future times}.

\subsect{Reducing a $(3 + 2n)$ DOF Polyped to a 3 DOF LLS}\label{sec:poly}

\begin{figure}[t]
\centering
\subfloat[polyped with $n=4$ legs,~leg~$k$~annotated]{%
\iftoggle{col2}{%
\def\svgwidth{0.65\columnwidth}%
\resizebox{!}{4cm}{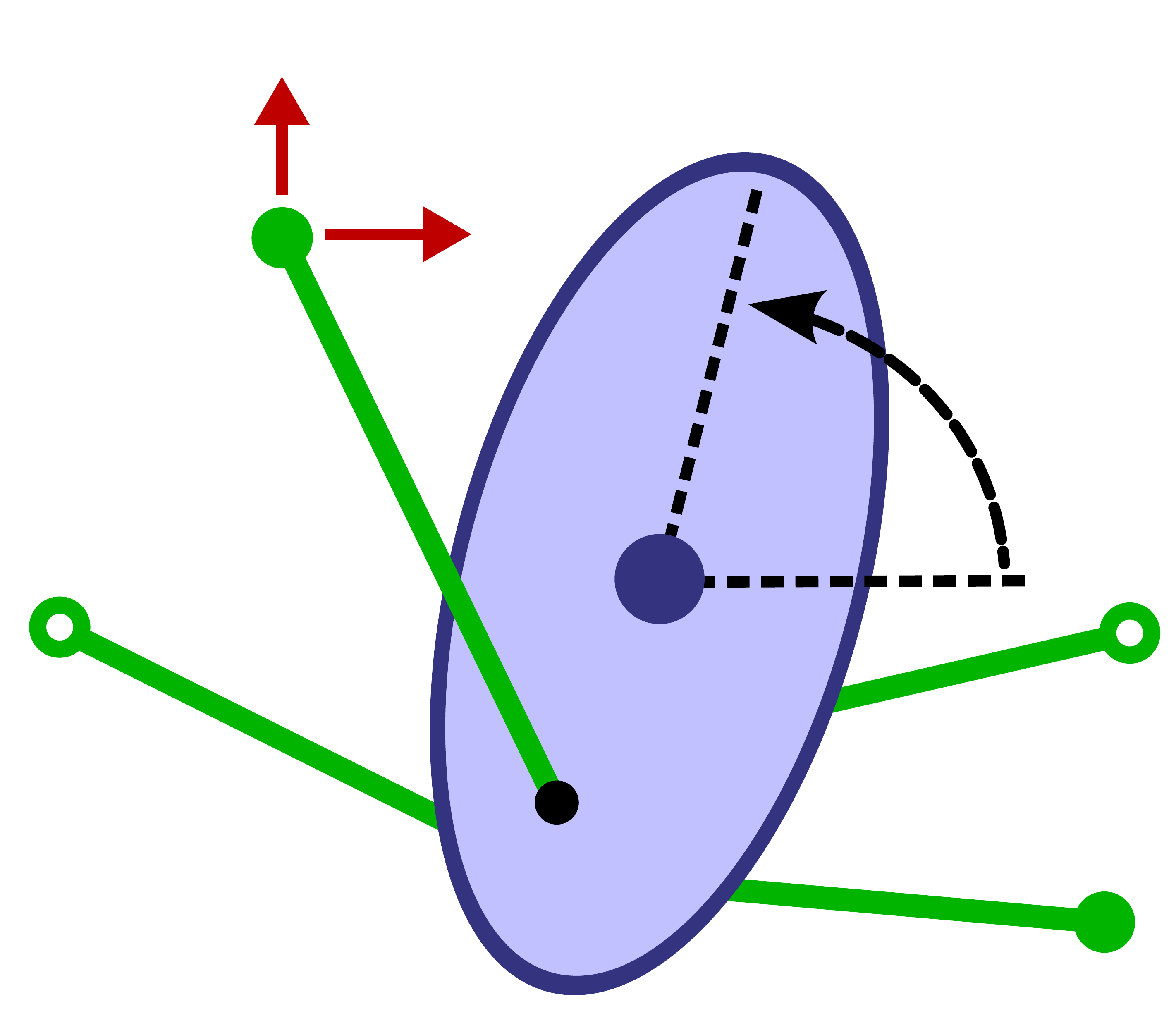}%
}{%
\def\svgwidth{0.45\columnwidth}%
\resizebox{.4\columnwidth}{!}{\input{poly.pdf_tex}}%
}%
\label{fig:poly}%
}%
\iftoggle{col2}{\hspace{.1cm}}{}%
\subfloat[Lateral Leg--Spring (LLS)~\cite{SchmittHolmes2000i}]{%
\iftoggle{col2}{%
\def\svgwidth{0.55\columnwidth}%
\resizebox{!}{4cm}{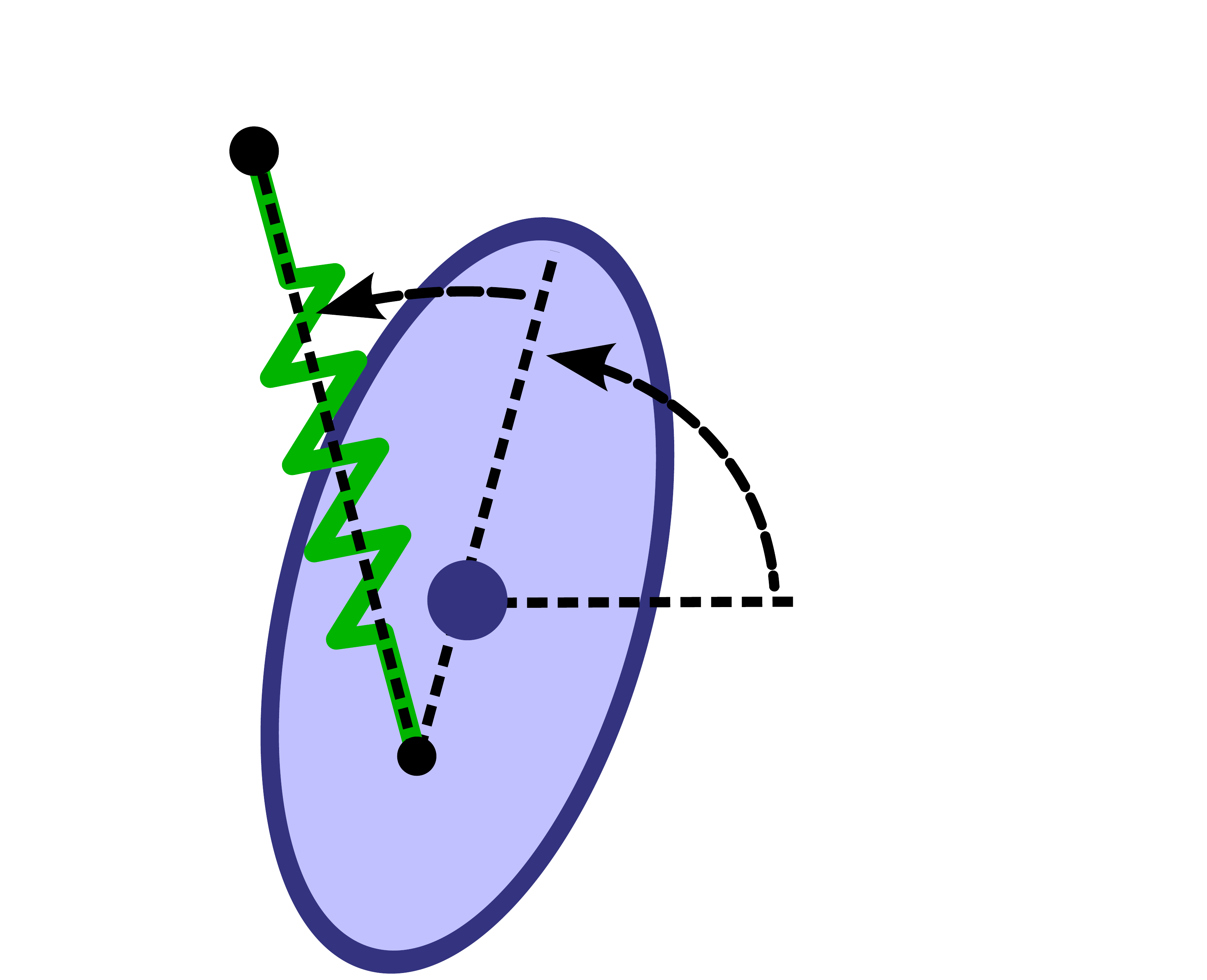}%
}{%
\def\svgwidth{0.5\columnwidth}%
\resizebox{.45\columnwidth}{!}{\input{lls.pdf_tex}}%
}%
\label{fig:lls}%
}%
\caption{
Lateral--plane models for locomotion described in Section~\ref{sec:poly}.
}\label{fig:legs}
\end{figure}
\revi{
A primary motivation for the present work is analysis of legged locomotion.
Several approaches have been proposed for embedding lower--dimensional dynamics in legged robot systems, notably \emph{hybrid zero dynamics}~\cite{WesterveltGrizzle2003} and \emph{active embedding}~\cite{AnkaraliSaranli2011}.
Complementing these engineering approaches and predating them, the \emph{templates and anchors hypotheses} (TAH)~\cite{FullKoditschek1999} conjectures that animal locomotion behaviors arise through reduction of the \emph{anchor} dynamics governing the nervous system and body~\cite{HolmesEtAl2006} to lower--dimensional \emph{template} dynamics that encode a specific behavior~\cite{SchmittHolmes2000i, GhigliazzaAltendorfer2003}.
One well--studied template is the Lateral Leg Spring (LLS)~\cite{SchmittHolmes2000i} model for sprawled posture running, which has been shown to match how cockroaches run and begin to recover from perturbations~\cite{JindrichFull2002}. 
Higher--dimensional neuromechanical variants of the model have been shown to reduce states associated with the nervous system~\cite{HolmesEtAl2006}.
In this section, we focus on reduction in the mechanical dynamics of limbs.
Specifically, we synthesize a state--feedback control law under which the underactuated lateral--plane polyped illustrated in Fig.~\ref{fig:poly} exactly reduces to the Lateral Leg--Spring (LLS)~\cite{SchmittHolmes2000i} model in Fig.~\ref{fig:lls}.
With $n$ limbs, the polyped possesses $(3+2n)$ degrees--of--freedom (DOF); the LLS has 3 DOF.
This example serves a dual purpose:
first, it demonstrates how our theoretical results can be applied to reduce an arbitrary number of DOF in a locomotion model; 
second, it suggests a mechanism that legged robot controllers could exploit to \emph{anchor} a desired \emph{template}. 

Before we proceed with describing the reduction procedure in detail, we give an overview of the approach and the connection with Theorem~\ref{thm:exact}.
We begin in Section~\ref{sec:poly:mdls} by describing the dynamics of the LLS template and polyped anchor.
Then in Section~\ref{sec:poly:embed} we construct a smooth state feedback law that ensures that trajectories of the polyped body exactly match those of the LLS; we accomplish this by simply ensuring the net \emph{wrench}~\cite{MurrayLi1994} comprised of generalized forces and torques acting on the polyped body matches that of the LLS for all time.
Subsequently, in Section~\ref{sec:poly:reduce} we modify the feedback law to further ensure the states associated with the polyped's limbs reduce after a single stride via Theorem~\ref{thm:exact}.
Finally, in Section~\ref{sec:poly:perturb} we discuss the effect of perturbations on the closed--loop reduced--order system.

\subsubsection{Dynamics of Lateral Leg--Spring (LLS) and $n$--leg polyped}\label{sec:poly:mdls}
The LLS is an energy--conserving lateral--plane model for locomotion comprised of a massless leg--spring with elastic potential $V$ affixed at hip position $h$ to an inertial body with two translational ($x,y$) and one rotational ($\theta$) DOF.
The system is initialized at the start of a stride by orienting the leg at a fixed angle $\beta$ with respect to the body at rest length $\ell$ and touching the foot down such that the leg will instantaneously contract.
The step ends once the leg extends to its rest length by touching the foot down on the opposite side of the body; subsequent steps are defined inductively.
In certain parameter regimes, the model possesses a periodic running gait~\cite{SchmittHolmes2000i}.

The underactuated hybrid control system illustrated in Fig.~\ref{fig:poly} extends neuromechanical models previously proposed to study multi--legged locomotion~\cite{HolmesEtAl2006, KukillayaProctor2009} by introducing masses into $n \ge 4$ feet connected by massless limbs affixed at hip locations $\set{h_k}_{k=1}^n$ on the inertial body.
We assume that each foot can attach or detach from the substrate at any time, and the transition from \emph{swing} to \emph{stance} entails a plastic impact that annihilates the kinetic energy in a foot.
We assume that each limb $k$ is fully--actuated; for simplicity we assume the inputs act along the Cartesian coordinates and do not saturate so that any $(\mu_k,\nu_k)\in\R^2$ is feasible at any limb configuration.
We let $q_0 = (x,y,\theta) \in Q_0 = \R^2\times S^1$ denote the position and orientation of the body, 
and for each $k\in\set{1,\dots,n}$ we let $q_k = (x_k,y_k) \in Q_k\in\R^2$ denote the position of the $k$--th foot.
The configuration space of the polyped is the $(n+1)$--fold product
$\prod_{k=0}^n Q_k$.
The $n$--leg polyped's dynamics thus have the form
\eqnn{\label{eqn:poly:dyn}
M \ddot{q}_0 = \sum_{k=1}^n (-\mu_k,-\nu_k,0) \Ad_{g_k},\ 
m_k \ddot{q}_k = (\mu_k,\nu_k)
}
where $M = \diag\paren{m,m,J}\in\R^{3\times 3}$ is the mass distribution of the body and $\Ad_{g_k}\in\R^{3\times 3}$ transforms a wrench applied at the $k$--th hip to an equivalent wrench applied at the body center--of--mass~\cite[\S5.1]{MurrayLi1994}.

\subsubsection{Embedding LLS in polyped}\label{sec:poly:embed}
For any subset $K\subset\set{1,\dots,n}$ of limbs, let 
\eqnn{\label{eqn:poly:wrench}
\sum_{k\in K} (-\mu_k,-\nu_k,0) \Ad_{g_k}\in T^*Q_0
}
denote the net wrench~\cite{MurrayLi1994} on the body resulting from actuating legs in $K$.
Then so long as no two hips are coincident \revii{and $\abs{K} \ge 2$},
any desired wrench 
may be imposed on the body by appropriate choice of inputs to the limbs in $K$ \revii{regardless of whether $K$ contains limbs in stance or swing}.
In the next section we describe a limb coordination procedure that ensures there will be a subset 
of stance limbs that can impose the LLS's wrench 
and cancel the reaction wrench from actuating other limbs 
at any time.

\subsubsection{Reducing polyped to LLS}\label{sec:poly:reduce}
We construct a smooth state feedback control law yielding a closed--loop Poincar\'{e} map $P_A:U_A\into\Sigma_A$ for the polyped 
that splits as $P_A:U_T\times U_N\into \Sigma_T\times\Sigma_N$ such that
\eqnn{\label{eqn:poly:TA}
P_A(z,\zeta) = (P_T(z), P_N(z))
}
where $P_T:U_T\into\Sigma_T$ is a Poincar\'{e} map for the LLS 
and $P_N:U_T\into \Sigma_N$ is a smooth map.
In the form~\eqref{eqn:poly:TA} it is clear that since $P_T$ is a diffeomorphism near the fixed point $\xi = P_T(\xi)$, all iterates of $P_A$ have constant rank equal to $\rank DP_T(\xi) = \dim\Sigma_T$ near $\xi$, and therefore Theorem~\ref{thm:exact} applies.

Partition the $n \ge 4$ limbs into $\swing\coprod\stance$, ensuring $\abs{\swing}$, $\abs{\stance}\ge 2$.  
Initialize at the beginning of a step at time $t$ with LLS and polyped body state $(q_0(t),\dot{q}_0(t))$ and polyped limb states $\set{(q_k(t),\dot{q}_k(t))}_{k=1}^n$ by attaching $\stance$ limbs and detaching $\swing$ limbs from the ground.
Note that the termination time $\tau$ for the LLS step depends smoothly on the initial condition $(q_0(t),\dot{q}_0(t))$.
For each $k\in\swing$ choosing constant inputs
\eqnn{\label{eqn:poly:ctrl}
\paren{\mu_k,\ \nu_k} = \frac{2}{\tau^2} \paren{(x(\tau),y(\tau)) + r(\theta(\tau)) \bar{q}_k - q_k(t) - \tau\dot{q}_k(t)}
}
ensures that the limb will reach a fixed location $\bar{q}_k$ in the body frame of reference at time $\tau$.
For each $k\in\stance$ choose inputs $(\mu_k,\nu_k)$ to cancel the reaction wrench from the swing limbs and impose the LLS acceleration on the polyped body.
At time $t+\tau$, exchange the $\stance$ and $\swing$ limb sets and proceed as with the previous step from the new initial condition.
After two steps, it is clear that the positions and velocities of the polyped's $n$ limbs are uniquely determined by the body initial condition $(q_0(t),\dot{q}_0(t))$. 
Therefore the polyped's Poincar\'{e} map has the form of~\eqref{eqn:poly:TA}, so Theorem~\ref{thm:exact} implies the polyped anchor reduces exactly to the LLS template after a single stride.


\subsubsection{Qualitative description of reduction}\label{sec:poly:qual}
The active embedding described in Section~\ref{sec:poly:embed} ensures the polyped body motion is always identical to that of the LLS, regardless of the state of the limbs.
The limb posture control in Section~\ref{sec:poly:reduce} guarantees the limb states are determined by the LLS body state after two steps, and furthermore synchronizes touchdown and liftoff events with those of the LLS.

\subsubsection{Effect of perturbations and parameter variations}\label{sec:poly:perturb}
The qualitative description in the preceding section makes it clear that, following a sufficiently small perturbation or parameter variation, the closed--loop polyped will continue to track and ultimately reduce to an LLS that experiences the corresponding disturbance.
Note that this conclusion requires that the polyped maintains the same control architecture exploited above to obtain the product decomposition in~\eqref{eqn:poly:TA}.
In particular, the controller must maintain observability of the full state and controllability of the limbs.
We study the effect of more general perturbations in the next section.
}

\subsect{Deadbeat Control of Rhythmic Hybrid Systems}\label{sec:deadbeat}

Generalizing the example from the previous section, we now consider a system wherein a finitely--parameterized control input updates when an execution passes through a distinguished subset of state space.
This form of control in rhythmic hybrid systems dates back (at least) to Raibert's hoppers~\cite{Raibert1986} and Koditschek's jugglers~\cite{BuehlerKoditschek1994}, and has recently received renewed interest~\revii{\cite{CarverCowan2009, Remy2011, PoulakakisGrizzle2009, AnkaraliSaranli2011, SeyfarthGeyer2003, SeipelHolmes2007, SreenathPark2013}}.
We model this with a hybrid system $H = (D, F, G, R)$ whose vector field and reset map depend on a control input that takes values in a smooth boundaryless 
manifold $\Theta$.
The value of the control input may be updated whenever an execution passes through the guard $G$, but it does not change in response to the continuous flow.
Suppose for some $\theta\in\Theta$ that $H$ possesses a periodic orbit $\gamma$,
let $P:U\times\Theta\into\Sigma$ be a \Pmap\, associated with $\gamma$ where $U\subset\Sigma\subset G$, and let $\set{\xi} = \gamma\cap\Sigma$.
\revi{
In this section we study deadbeat control of the discrete--time nonlinear control system
\eqnn{\label{eqn:dt:nl}
x_{i+1} = P(x_i, \theta_i)
}
and the discrete--time linear control system obtained by linearizing $P$ about the fixed point $\xi = P(\xi,\theta)$,
\eqnn{\label{eqn:dt:l}
\dx_{i+1} = D_x P(\xi,\theta) \dx_i + D_\theta P(\xi,\theta) \dtheta_i.
}
The control architecture we present is well--known for linear and nonlinear maps arising in locomotion~\cite{CarverCowan2009}; 
the novelty of this section lies in the connection to \emph{exact} and \emph{approximate} reduction via Theorems~\ref{thm:exact} and~\ref{thm:approx}.

\subsubsection{Exact reduction over one cycle}\label{sec:onecycle}
As studied in~\cite{CarverCowan2009}, an application of the Implicit Function Theorem~\cite[Theorem~7.8]{Lee2002} shows that if $\rk D_\theta P(\xi,\theta) = \dim\Sigma$ then there exists a neighborhood $V\subset U$ of $\xi$ and a smooth feedback law $\psi:V\into\Theta$ such that for all $x\in V$ we have $P(x,\psi(x)) = \xi$, i.e. $\psi$ is a \emph{deadbeat} control law for~\eqref{eqn:dt:nl}.
Since $\psi$ is smooth, the closed--loop \Pmap\, $P_\psi:V\into\Sigma$ defined by $P_\psi(x) = P(x,\psi(x))$ satisfies the hypotheses of Theorem~\ref{thm:exact} (Exact Reduction) with rank $r = 0$, so the invariant subsystem yielded by the Theorem is simply the periodic orbit $\gamma$.

In practice it may be desirable to reduce fewer than $\dim\Sigma$ coordinates.
If there exists a smooth function $h:\Sigma\into\R^d$ that satisfies $h\circ P(\xi,\theta) = 0$ and $\rank D_\theta h\circ P(\xi,\theta) = d$, then the preceding construction yields a closed--loop system that reduces via Theorem~\ref{thm:exact} to the embedded $d$--dimensional submanifold $h^{-1}(0)$ near $\xi$.

\subsubsection{Exact reduction over multiple cycles}\label{sec:manycycles}
If $\rank D_\theta P(\xi,\theta) < \dim\Sigma$, as noted in~\cite{CarverCowan2009} it may be possible to construct a deadbeat control law by applying inputs over multiple cycles.
Specifically, let $P_0 = P$ and for each $\ell\in\N$ define $P_\ell:U_\ell\times\Theta^\ell\into\Sigma$ by 
\eqnn{\label{eqn:iterate}
P_{\ell}(x,(\theta_1,\dots,\theta_{\ell})) = P(P_{\ell-1}(x,(\theta_1,\dots,\theta_{\ell-1})), \theta_{\ell})
}
for all $(x,(\theta_1,\dots,\theta_{\ell}))\in U_\ell\times\Theta^\ell$
where $U_\ell\subset U$ is a neighborhood of $\xi$ sufficiently small to ensure~\eqref{eqn:iterate} is well--defined.
Then if there exists $k\in\N$ such that
\eqnn{\label{eqn:ctrlability}
\rank D_{(\theta_1,\dots,\theta_k)} P_k(\xi,(\theta,\dots,\theta)) = \dim\Sigma,
}
the construction from the previous paragraph yields a smooth $k$--step feedback law $\psi_k:V_k\into\Theta^k$ such that the closed--loop hybrid system reduces via Theorem~\ref{thm:exact} to the periodic orbit $\gamma$ after $k$ cycles.
We conclude this section by noting that~\cite{CarverCowan2009} contains an example that performs exact reduction after two cycles, and for which reduction in fewer cycles is impossible.

\subsubsection{Approximate reduction}\label{sec:controllability}
Since~\eqref{eqn:ctrlability} is equivalent to controllability~\cite[Chapter~8d.5]{CallierDesoer1991} of the linear control system~\eqref{eqn:dt:l},
it is worthwhile to consider the linear control problem. 
Any stabilizable subspace $S$~\cite[Chapter~8d.7]{CallierDesoer1991} of~\eqref{eqn:dt:l} can be rendered attracting in a finite number of steps $k\in\N$ with linear state feedback $\dtheta_i = \Psi \dx_i$ where $\Psi$ is a fixed matrix~\cite{OReilly1981}.
Applying this linear feedback law to the nonlinear system~\eqref{eqn:dt:nl} yields a closed--loop Poincar\'{e} map $P_\Psi$ such that the rangespace of the $k$--th iterate of its linearization $D_x P^k_\Psi(\xi)$ is contained in $S$. 
Therefore Theorem~\ref{thm:approx} (Approximate Reduction) yields an invariant hybrid subsystem, tangent to $S$ on $\Sigma$, that attracts nearby trajectories superexponentially. 
Thus, although feedback laws for the nonlinear control system~\eqref{eqn:dt:nl} constructed above can be computed using the procedure described in~\cite{CarverCowan2009} to achieve exact reduction to the target subsystem, if approximate reduction suffices then one may simply apply the linear deadbeat controller computed for~\eqref{eqn:dt:l}.
}
\subsubsection{Structural stability of deadbeat control}\label{sec:structstab}
Suppose the preceding development is applied to a model that differs from that used to construct the feedback law $\psi\in C^\infty(V,\Theta)$.
We study the \emph{structural stability}~\cite[Section~1.7]{GuckenheimerHolmes1983} of attracting invariant sets arising in this class of systems by applying the Theorems of Section~\ref{sec:hds}.
If the models differ by a small smooth deformation (as would occur if there was a small perturbation in model parameters), one interpretation of this change is that some $\td{\psi}\in B_\veps(\psi)\subset C^\infty(V,\Theta)$ is applied to the model for which $\psi$ is deadbeat, where $\veps > 0$ bounds the error.
For all $\veps > 0$ sufficiently small, $\td{\psi}$ yields a perturbed closed--loop \Pmap\, $\td{P}:V\into\Sigma$ possessing a unique fixed point $\td{\xi}\in V$, and 
$\td{\xi}$ is an exponentially stable fixed point of the perturbed system. 

We conclude by noting that it is possible for the structure of the hybrid dynamics to constrain the achievable perturbations.
For instance, if one domain of the hybrid system has lower dimension than that in which the \Pmap\, is constructed, then zero is always a Floquet multiplier regardless of the applied feedback; in this case Theorem~\ref{thm:approx} (Approximate Reduction) implies the existence of a proper submanifold of the Poincar\'{e} section $\Sigma$ to which trajectories contract superexponentially in the presence of any (sufficiently small) smooth perturbation to the closed--loop dynamics.

\subsect{Hybrid Floquet Coordinates}\label{sec:floq}

When a hybrid system reduces to a smooth dynamical system near a periodic orbit
via Theorem~\ref{thm:exact}~(Exact Reduction), we can generalize the \emph{Floquet normal form}~\cite{floquet1883edl, Guckenheimer1975, Revzen2009, RevzenGuckenheimer2011} using Theorem~\ref{thm:smooth}~(Smoothing).
\revi{
Broadly, this demonstrates how the Theorems of Section~\ref{sec:hds} can be applied to generalize constructions from classical dynamical systems theory to the hybrid setting.
More concretely, this provides a theoretical framework that justifies application of the empirical approach developed in~\cite{Revzen2009, RevzenGuckenheimer2011} to estimate low--dimensional invariant dynamics in data collected from physical locomotors.
}

Consider a hybrid dynamical system $H = (D,F,G,R)$ with $\tau$--periodic orbit $\gamma$ that satisfies the hypotheses of Theorem~\ref{thm:exact}.
Let $M\subset D$ be the $(r+1)$--dimensional invariant hybrid subsystem yielded by the Theorem, and $W\subset D$ a hybrid open set containing $\gamma$ that contracts to $M$ in finite time.
Let $(\td{M},\td{F})$ denote the smooth dynamical system obtained by applying Theorem~\ref{thm:smooth}.
Under a genericity condition\footnote{Either the periodic orbit is exponentially stable or it is \emph{hyperbolic} and the associated \emph{Floquet multipliers} do not satisfy any \emph{Diophantine equation}~\cite[Chapter~3.3]{GuckenheimerHolmes1983}.} there exists a neighborhood $U\subset\td{M}$ of $\gamma$ and a smooth chart $\vphi:U\into\R^r\times S^1$ such that the coordinate representation of the vector field has the form
\eqnn{\label{eqn:floq}
D\vphi\circ\td{F}\circ D\vphi^{-1}(z,\theta) = \mat{c}{\dot{z} \\ \dot{\theta}} = \mat{c}{A(\theta) z \\ 2\pi/\tau}
}
where $z\in\R^r$ and $\theta\in S^1$.
In these coordinates, each $\theta\in S^1$ determines an embedded submanifold $\td{N}_\theta = \R^r\times\set{\theta}\subset\R^r\times S^1$ that is mapped to itself after flowing forward in time by $\tau$; for this reason, the submanifolds $\td{N}_\theta$ are referred to as \emph{isochrons}~\cite{Guckenheimer1975}.
Each $x\in \td{N}_\theta$ may be assigned the \emph{phase} $\theta\in S^1$; if $\gamma$ is stable, then as $t\goesto\infty$ the trajectory initialized at $x$ will asymptotically converge to the trajectory initialized at $(0,\theta)$.

The isochrons may be pulled back to any precompact hybrid open set $V\subset W$ containing $\gamma$ in the original hybrid system as follows.
The proof of Theorem~\ref{thm:exact} implies there exists a finite time $t < \infty$ such that every execution initialized in $V$ is defined over the time interval $[0,t]$ and reaches $M$ before time $t$; without loss of generality, we take this time to be a multiple $k\tau$ of the period of $\gamma$ for some $k\in\N$.
Let $\psi:V\into \td{M}$ denote the map that flows an initial condition $x\in V$ forward by $t$ time units and then applies the quotient projection $\pi:M\into\td{M}$ obtained from Theorem~\ref{thm:smooth} to yield the point $\psi(x)\in \td{M}$.
Then the constructions in the proof of Theorem~\ref{thm:exact} imply that $\psi$ is a smooth map in the sense defined in Section~\ref{sec:hdg}, i.e. it is continuous and $\psi|_{V\cap D_j}$ is smooth for each $j\in J$.
Now for any $\theta\in S^1$ the set $N_\theta = \psi^{-1}(U)$ is mapped into $\td{N}_\theta$ after $k\tau$ units of time; we thus refer to $N_\theta\subset D$ as a \emph{hybrid isochron}.
We conclude by noting that $N_\theta$ will generally not be a smooth (hybrid) submanifold.

\sect{Discussion}\label{sec:disc}
Generically for an exponentially stable periodic orbit in a hybrid dynamical system, nearby trajectories contract superexponentially to a subsystem containing the orbit.
Under a non--degeneracy condition on the rank of any \Pmap\, associated with the orbit, this contraction occurs in finite time regardless of the stability of the orbit.
Hybrid transitions may be removed from the resulting subsystem, yielding an equivalent smooth dynamical system.
Thus the dynamics near stable hybrid periodic orbits are generally obtained by extending the behavior of a smooth system in transverse coordinates that decay superexponentially.
Although the applications presented in Section~\ref{sec:ex} focused on terrestrial locomotion~\cite{HolmesEtAl2006}, we emphasize that the results in Section~\ref{sec:hds} do not depend on the phenomenology of the physical system under investigation, and are hence equally suited to study rhythmic hybrid control systems appearing in robotic manipulation~\cite{BuehlerKoditschek1994}, biochemistry~\cite{GlassPasternack1978}, and electrical systems~\cite{HiskensReddy2007}.

In addition to providing a canonical form for the dynamics near hybrid periodic orbits, the results of this paper suggest a mechanism by which a many--legged locomotor or a multi--fingered manipulator may collapse a large number of mechanical degrees--of--freedom to produce a low--dimensional coordinated motion.
This provides a link between disparate lines of research: formal analysis of hybrid periodic orbits; design of robots for rhythmic locomotion and manipulation tasks; and scientific probing of neuromechanical control architectures in humans and animals.
\revi{
Our theoretical results show that hybrid models of rhythmic phenomena generically reduce dimensionality, and our applications demonstrate that this reduction may be deliberately designed into an engineered system.
We furthermore speculate that evolution may have exploited this reduction in developing its spectacularly dexterous agents.
}


\subsection*{Support}
This research was supported in part by: an NSF Graduate Research Fellowship to S. A. Burden;
ARO Young Investigator Award \#61770 to S. Revzen;
and Army Research Laboratory Cooperative Agreements W911NF--08--2--0004 and W911NF--10--2--0016.
The views and conclusions contained in this document are those of the authors and should not be interpreted as representing the official policies, either expressed or implied, of the Army Research Laboratory or the U.S. Government.  
The U.S. Government is authorized to reproduce and distribute for Government purposes notwithstanding any copyright notation herein.

\subsection*{Acknowledgments}
We thank John Hauser for correcting an error in an earlier version of this manuscript, and
Sam Coogan,
John Guckenheimer, 
Ram Vasudevan,
and the three anonymous reviewers
for their invaluable feedback.

\appendices

\sect{Smooth Dynamical Systems}\label{sec:sds}

We constructed hybrid systems using switching maps defined on boundaries of smooth dynamical systems.  
The behavior of such systems can be studied by alternately applying flows and maps,
thus in this section we collect results that provide canonical forms for the behavior of flows and maps near periodic orbits and fixed points, respectively.
The first develops a canonical form for the flow to a section in a continuous--time system.
The second provides a technique to smoothly attach continuous--time systems along their boundaries.
The third and fourth establish a canonical form for submanifolds that are invariant and approximately invariant (respectively) near fixed points in discrete--time dynamical systems;
the fifth provides an estimate of the error in the invariance approximation.

\subsect{Continuous--Time Dynamical Systems}\label{sec:ct}

\defn{A \emph{continuous--time dynamical system} is a pair $(M,F)$ where: 
\begin{itemize}
\item[$M$] is a smooth manifold with boundary $\bd M$;
\item[$F$] is a smooth vector field on $M$, i.e. $F\in\e{T}(M)$.
\end{itemize}}

\subsubsection{Time--to--Impact}\label{sec:tti}
When a trajectory passes transversely through an embedded submanifold,
the time required for nearby trajectories to pass through the manifold depends smoothly on the initial condition~\cite[Chapter 11.2]{HirschSmale1974}. 
This provides the prototype used in the proofs of Theorems~\ref{thm:exact} and~\ref{thm:approx} for the dynamics near the portion of a hybrid periodic orbit in one domain of a hybrid system. 

\lem{\label{lem:sec}
Let $(M,F)$ be a smooth dynamical system, $\phi:\e{F}\into M$ the maximal flow associated with $F$, and $G\subset M$ a smooth codimension--1 embedded submanifold.
If there exists $x\in M$ and $t\in\e{F}^x$ such that $\phi(t,x)\in G$ and $F\paren{\phi(t,x)}\not\in T_{\revii{\phi(t,x)}}G$, then there is a neighborhood $U\subset M$ containing $x$ and a smooth map $\sigma:U\into\R$ so that $\sigma(x) = t$ and $\phi(\sigma(y),y)\in G$ for all $y\in U$;  $\sigma$ is called the \emph{time--to--impact} map.
}
\iftoggle{sds}{
\pf{
Near $\phi(t,x)$, $G$ is the zero section of a constant--rank map $h:M\into\R$ where $Dh\paren{\phi(t,x)}\ne 0$.  Define $g:\e{F}\into\R$ by $g(s,y) = (h\circ\phi)(s,y)$.  Then since $F$ is transverse to $G$ at $\phi(t,x)$, 
\[ \vfof{g}{t}(t,x) = Dh\paren{F\paren{\phi(t,x)}} \ne 0. \]
By the Implicit Function Theorem \see{Theorem 7.8 in~\cite{Lee2002}}, there exists a neighborhood $U$ of $x$ and a smooth map $\sigma:U\into\R$ so that $\sigma(x) = t$ and $g(\sigma(y),y) = 0$ for all $y\in U$, i.e. $\phi(\sigma(y),y) \in G$.
}
}

\rem{This lemma is applicable when $G\subset\bd M$.}

\subsubsection{Smoothing Flows}\label{sec:sds:smooth}
Two continuous--time dynamical systems can be smoothly attached to one another along their boundaries to obtain a new continuous--time system~\cite[Theorem~8.2.1]{Hirsch1976}.
Distinct hybrid domains were attached to one another using this construction in Section~\ref{sec:hds}.

\lem{\label{lem:smooth}
Suppose $(M_1,F_1),(M_2,F_2)$ are $n$--dimensional continuous--time dynamical systems, 
there exists a diffeomorphism
$R:\bd M_1\into \bd M_2$, 
$F_1$ is outward--pointing along $\bd M_1$, and $F_2$ is inward--pointing along $\bd M_2$.  
Then the topological quotient 
\[
\td{M} = \frac{M_1\coprod M_2}{\bd M_1\overset{R}{\sim}\bd M_2}
\]
can be endowed with the structure of a smooth manifold such that for $j\in\set{1,2}$:
\begin{enumerate}
\item the quotient projections $\pi_j:M_j\into \td{M}$ are smooth embeddings; and
\item there is a smooth vector field $\td{F}\in\e{T}(\td{M})$ that restricts to $D\pi_j(F_j)$ on $\pi(M_j)\subset\td{M}$.
\end{enumerate}
}

\iftoggle{sds}{
\pf{
Let $\phi_j:\e{F}_j\into M_j$ be the maximal flow associated with $F_j$ on $M_j$.
Then there is a neighborhood $\td{U}_j\subset\e{F}_j$ of $\set{0}\times M_j$ on which the flow is defined, and with $U_j = \td{U}_j\cap (\R\times\bd M_j)$, transversality of $F_j$ along $\bd M_j$ implies $\phi_j:U_j\into M_j$ is an embedding which is the identity on $\set{0}\times\bd M_j$.  
Since $F_1$ is outward--pointing and $F_2$ is inward--pointing, the neighborhoods are one--sided and without loss of generality we may assume for $j = 1,2$ that there exist continuous positive functions $\delta_j:\bd M_j\into[0,\infty)$ such that $U_1 = \set{(-\delta_1(x),0] : x\in\bd M_1}$ and $U_2 = \set{[0,\delta_2(x)) : x\in\bd M_2}$.
Therefore $U = \frac{U_1\coprod U_2}{\bd M_1\simeq\bd M_2}$ inherits a smooth structure from its product structure, i.e. the fibers $U^x = (-\delta_1(x),\delta_2(\vphi(x)))\times\set{x}$ are smooth curves for $x\in\bd M_1$ and both $\set{0}\times\bd M_1\inc U$ and $\set{0}\times \vphi^{-1}(\bd M_2)\inc U$ are smooth embeddings; let $\vphi:U\into\R^n$ denote the chart.
Note in addition that by construction the constant vector field $\vf{t}\in\e{T}(U_j)$ pushes forward to $F_j|_{\phi_j(U_j)\cap M_j}\in\e{T}\paren{\phi_j(U_j)\cap M_j}$, since $(D\phi_j)\vf{t} = F_j$.

We construct the smooth structure on $\td{M} = \frac{M_1\coprod M_2}{\bd M_1 \simeq\bd M_2}$ by covering $M$ with interior charts from the $M_j$'s together with $U$.  
Note that since $\phi_j|_{U_j}:U_j\into M_j$ is a smooth embedding, the interior charts on $M_j$ are smoothly compatible with the product chart on $U$, and the natural quotient projections $\pi_j:M_j\inc \td{M}$ are smooth embeddings.  
Finally, $(D\pi_1) F_1 = (D\pi_2) F_2$ along $\pi_1(\bd M_1)$ by construction,
whence the vector field $\td{F}\in\e{T}(\td{M})$ which restricts to $F_j$ on $M_j$, $j=1,2$, is well--defined and smooth.
}
}

\rem{The smooth structure constructed in Lemma~\ref{lem:smooth} is unique up to diffeomorphism~\cite[Theorem~2.1 in Chapter~8]{Hirsch1976}.
}

\subsect{Discrete--time Dynamical Systems}\label{sec:dt}
\defn{A \emph{discrete--time dynamical system} is a pair $(\Sigma,P)$ where: 
\begin{itemize}
\item[$\Sigma$] is a smooth manifold without boundary;
\item[$P$] is a smooth endomorphism of $\Sigma$, i.e. $P:\Sigma\into\Sigma$.
\end{itemize}}

In studying hybrid dynamical systems, we encounter smooth maps $P:\Sigma\into \Sigma$ that are noninvertible.
Viewing iteration of $P$ as determining a discrete--time dynamical system, we wish to study the behavior of these iterates near a fixed point $\xi = P(\xi)$.
Note that if $P$ has constant rank equal to $k < n = \dim \Sigma$, then its image $P(\Sigma)\subset \Sigma$ is an embedded $k$--dimensional submanifold near $\xi$ by the Rank Theorem~\cite[Theorem~7.13]{Lee2002}.
With an eye toward model reduction, one might hope that the composition $(P\circ P):\Sigma\into P(\Sigma)$ is also constant--rank, but this is not generally true\footnote{Consider the map $P:\R^2\into\R^2$ defined by $P(x,y) = (x^2,x)$.}.

In this section we provide three results that introduce regularity into iterates of a noninvertible map $P:\Sigma\into \Sigma$ on an $n$--dimensional manifold $\Sigma$ near a fixed point $P(\xi) = \xi$.
If the rank of $DP$ is strictly bounded above by $m\in\N$ and 
if $P^m$, the $m$--th iterate of $P$, has constant rank equal to $r\in\N$ near the fixed point $\xi$, then $P$ reduces to a diffeomorphism over an $r$--dimensional invariant submanifold after $m$ iterations; this result is given in Section~\ref{sec:sds:exact}.
Even if $DP^m$ is not constant rank, as long as $\xi$ is exponentially stable then $P$ can be approximated by a diffeomorphism on a submanifold whose dimension equals $\rank DP^m(\xi)$; this is the subject of Section~\ref{sec:sds:approx}.
A bound on the error in this approximation is provided in Section~\ref{sec:superstable}.

\subsubsection{Exact Reduction}\label{sec:sds:exact}
If the rank of $P:\Sigma\into\Sigma$ is strictly bounded above by $m\in\N$
and the derivative of the $m$--th iterate of $P$ has constant rank near a fixed point, then the range of $P$ is locally an embedded submanifold, and $P$ restricts to a diffeomorphism over that submanifold.
This originally appeared without proof as Lemma~3 in~\cite{BurdenRevzen2011}.

\lem{\label{lem:exact}
Let $(\Sigma,P)$ be an $n$--dimensional discrete--time dynamical system with $P(\xi) = \xi$ for some $\xi\in \Sigma$. 
Suppose the rank of $P$ is \revi{strictly} bounded above by $m\in\N$ and 
\revi{there exists a neighborhood $W\subset\Sigma$ of $\xi$ such that $\rank DP^m(x) = r$ for all $x\in W$}. 
Then there is a neighborhood $V\subset \Sigma$ containing $\xi$ such that $P^m(V)$ is an $r$--dimensional embedded submanifold near $\xi$ and there is a neighborhood  $U\subset P^m(V)$ containing $\xi$ that $P$ maps diffeomorphically onto $P(U)\subset P^m(V)$.
}
\noindent
In the proof of Lemma~\ref{lem:exact}, we make use of a fact from linear algebra obtained by passing to the Jordan canonical form.
\revi{
\prop{\label{prop:la1}
If $A\in\R^{n\times n}$ and $\rk A < m$, then $\rk(A^{2m}) = \rk(A^m)$.
}
}

\pf{{\it (of Lemma~\ref{lem:exact})}
By the Rank Theorem~\cite[Theorem~7.13]{Lee2002}, there is a neighborhood $V\subset \Sigma$ of $\xi$ for which $S = P^m(V)$ is an $r$--dimensional embedded submanifold
and by Proposition~\ref{prop:la1} we have
\eqnal{
  \rk\paren{DP^m|_S}(\xi) & = \rk D(P^m\circ P^m)(\xi) 
  \\
  & = \rk DP^m(\xi).
}
Therefore $DP^m|_S:T_\xi S\into T_\xi S$ is a bijection, so by the Inverse Function Theorem~\cite[Theorem~7.10]{Lee2002}, there is a neighborhood $W\subset S$ containing $\xi$ so that $P^m(W)\subset S$ and $P^m|_{W}:W\into P^m(W)$ is a diffeomorphism.

We now show that $W$ is invariant under $P$ in a neighborhood of $\xi$.
By continuity of $P$, there is a neighborhood $L\subset V$ containing $\xi$ for which $P(L)\subset V$ and $P^m(L)\subset W$. 
The set 
$U = P^m(L)$
is a neighborhood of $\xi$ in $S$.
Further, we have 
\eqnal{
  P(U) & = P\circ P^m(L) = P^m\circ P(L)\subset S.
}
The restriction $P^m|_U:U\into P^m(U)$ is a diffeomorphism since $U \subset W$, whence $P|_U$ is a diffeomorphism onto its image $P(U)\subset S$.
}

\subsubsection{Approximate Reduction}\label{sec:sds:approx}
Now suppose that iterates of $P$ are not constant rank but $\xi = P(\xi)$ is exponentially stable, meaning that the \emph{spectral radius} $\specr{DP(\xi)} = \max\set{\abs{\lambda} : \lambda\in\spec DP(\xi)}$ satisfies $\specr{DP(\xi)} < 1$.
We show that $P$ may be approximated by a diffeomorphism defined on a submanifold whose dimension equals the number of non--zero eigenvalues of $DP(\xi)$.
The technical result we desire was originally established by Hartman~\cite{Hartman1960a}\footnote{The statement in~\cite{Hartman1960a} only considered invertible contractions.
However, as noted in~\cite{AulbachGaray1994}, the proof in~\cite{Hartman1960a} of the result we require does not make use of invertibility and the conclusion is still valid if zero is an eigenvalue of the linearization.
For details we refer to \cite{Abbaci2004}. 
}.
We apply Hartman's Theorem to construct a $C^1$ change--of--coordinates that exactly linearizes all eigendirections corresponding to non--zero eigenvalues of $DP(\xi)$.

\lem{\label{lem:approx}
Let $(\Sigma,P)$ be an $n$--dimensional discrete--time dynamical system.
Suppose $\xi = P(\xi)$ is an exponentially stable fixed point and let 
\revi{$r$ be the number of non--zero eigenvalues of $DP(\xi)$}%
.
Then there is a neighborhood $U\subset \Sigma$ of $\xi$ and a $C^1$ diffeomorphism $\vphi:U\into\R^n$ such that $\vphi(\xi) = 0$ and the coordinate representation $\td{P} = \vphi\circ P\circ\vphi^{-1}$ of $P$ has the form
\eqnal{
\td{P}(z,\zeta) = \paren{Az,\ N(z,\zeta)}
}
where $z\in\R^r$, $\zeta\in\R^{n-r}$, 
$A\in\R^{r\times r}$ is invertible, 
$N:\vphi(U)\into\R^{n-r}$ is $C^1$, 
$N(0,0) = 0$,
and $D_\zeta N(0,0)$ is nilpotent.
}

\iftoggle{sds}{
\pf{
Let $(U_0,\vphi_0)$ be a smooth chart for $\Sigma$ with $\xi\in U_0$ and $\vphi_0(\xi) = 0$.
We begin by verifying that the hypotheses of Theorem~\ref{thm:hartman} \revii{from Appendix~\ref{app:c1}} are satisfied for the map $P_0:\vphi_0(U_0)\into\R^n$ defined by $P_0 = \vphi_0\circ P\circ\vphi_0^{-1}$.
Let $\lambda\in\spec DP_0(0)$ be the eigenvalue with largest magnitude, and
$\ell\in\N$ its algebraic multiplicity.
Applying the linear change--of--coordinates that puts $DP_0(0)$ into Jordan canonical form, we assume 
\eqn{
DP_0(0) &= \mat{ccc}{A & 0 \\ 0 & B}
}
where 
$B\in\R^{\ell\times\ell}$ and $\spec B = \set{\lambda}$.
Now in the notation of Theorem~\ref{thm:hartman}, 
\eqn{
P_0(x,y) = \paren{Ax + X(x,y), By + Y(x,y)}
}
where $x\in\R^{n-\ell}$, $y\in\R^\ell$, and $X$, $Y$ are smooth and $X(0,0) = 0$, $Y(0,0) = 0$; note that $m = 0$ (there is no $z$ coordinate) at this step.
Because $X$ and $Y$ are smooth on the neighborhood $U_0$ of the origin, their derivatives are uniformly Lipschitz and H\"{o}lder continuous on a precompact open subset of $U_0$.

Theorem~\ref{thm:hartman} implies there exists a neighborhood $U_1\subset\R^n$ of the origin and a $C^1$ diffeomorphism $\vphi_1:U_1\into\R^n$ for which the map $P_1:\vphi_1(U_1)\into\R^n$ defined by $P_1 = \vphi_1\circ P_0\circ\vphi_1^{-1}$ has the form (after reversing the order of the coordinates)
\eqnal{
P_1(z_1,\zeta_1) = \paren{A_1z_1,\ N_1(z_1,\zeta_1)}
}
where $z_1\in\R^{r_1}$, $r_1 > 0$, $\zeta_1\in\R^{n-r_1}$ and $A_1\in\R^{r_1\times r_1}$ is invertible.
Observe that the map $P_1$ satisfies the hypotheses of Theorem~\ref{thm:hartman}.
Therefore we may inductively apply the Theorem to construct a sequence of coordinate charts $\set{(U_k,\vphi_k)}_{k=1}^K$ and corresponding maps $\set{P_k}_{k=1}^K$ such that for all $k\in\set{1,\dots,K}$
\eqnal{
P_k(z_k,\zeta_k) = \paren{A_kz_k,\ N_k(z_k,\zeta_k)}
}
where $z_k\in\R^{r_k}$, $\zeta_k\in\R^{n-r_k}$, $A_k\in\R^{r_k\times r_k}$ is invertible,
and $r_k > r_{k-1}$ (note that $r_0 = 0$).
The sequence terminates at a finite $K < \infty$ with $r_K = r = \rk DP^n(\xi)$.
Therefore in the $C^1$ chart $(U,\vphi)$ given by $\vphi = \vphi_K\circ\cdots\circ\vphi_0$ and $U = \vphi^{-1}(\R^n)$, the coordinate representation $\td{P} = \vphi\circ P\circ\vphi^{-1}$ of $P$ has the form
\eqnal{
\td{P}(z,\zeta) = \paren{Az,\ N(z,\zeta)}
}
where $z\in\R^r$, $\zeta\in\R^{n-r}$ and $A\in\R^{r\times r}$ is invertible. 
Since $A$ is invertible and $\rk D\td{P}^n(\xi) = r$, $D_\zeta N(0,0)$ is nilpotent.
}
}

\subsubsection{Superstability}\label{sec:superstable}

Finally, we recall that if all eigenvalues of the linearization of a map at a fixed point are zero---a so--called ``superstable'' fixed point~\cite{WendelAmes2012}---then the map contracts superexponentially;\iftoggle{sds}{\footnote{The map need not be nilpotent simply because its linearization is; consider the map $P:\R\into\R$ defined by $P(x) = x^2$.}}{}
this is a straightforward consequence of Ostrowski's Theorem~\cite[8.1.7]{Ortega1990}.

\lem{\label{lem:superstable}
Let $P:\R^n\into\R^n$ be a $C^1$ map with $P(0) = 0$, $\spec DP(0) = \set{0}$.
Then for every $\veps > 0$ and norm $\norm{\cdot}:\R^n\into\R$ there exists $\delta,C > 0$ such that
\eqnal{
\forall x\in B_\delta(0),k\in\N : \norm{P^k(x)} \le C\veps^k\norm{x}.
}
}
\iftoggle{sds}{
\noindent
The proof of Lemma~\ref{lem:superstable} relies on the following elementary fact regarding induced norms~\cite[1.3.6]{Ortega1990}.

\begin{proposition}[1.3.6 in~\cite{Ortega1990}]\label{prop:norm}
Given $\veps > 0$ and $A\in\R^{n\times n}$, there exists a norm $\norm{\cdot}:\R^n\into\R$ such that
$\norm{A}_i \le \specr{A} + \veps$,
where $\inorm{\cdot}:\R^{n\times n}\into\R$ is the operator norm induced by $\norm{\cdot}$ and $\rho(A)$ is the spectral radius of $A$.
\end{proposition}

\pf{(of Lemma~\ref{lem:superstable})
Given $\veps > 0$, choose the norm $\norm{\cdot}:\R^n\into\R$ obtained by applying Proposition~\ref{prop:norm} to $DP(0)$ so that $\inorm{DP(0)}\le\frac{1}{2}\veps$.
Since $DP$ is continuous, there exists a $\delta > 0$ such that
\eqnal{
\forall x\in B_\delta(0) : \inorm{DP(x) - DP(0)} < \frac{1}{2}\veps.
}
Whence we find for $\norm{x} < \delta$ that
\eqnal{
\inorm{DP(x)} &= \inorm{DP(x) - DP(0) + DP(0)} \\
&\le \inorm{DP(x) - DP(0)} + \inorm{DP(0)} \le \veps. \\
}
Combined with 8.1.4 in~\cite{Ortega1990} (a generalization of the Mean Value Theorem to vector--valued functions), we find for all $x\in B_\delta(0)$,
\eqnal{
\norm{P(x)} &\le \sup_{s\in[0,1]}\inorm{DP(sx)}\norm{x} \le \veps\norm{x}.
}
Iterating, for all $k\in\N$ and $\norm{x}<\delta$ we have $\norm{P^k(x)}\le\veps^k\norm{x}$.
Since all norms on finite--dimensional vector spaces are equivalent, the desired result follows immediately. 
}
}

\rem{
Let $(\Sigma,P)$ be an $n$--dimensional discrete--time dynamical system that satisfies the hypotheses of Lemma~\ref{lem:approx} near $\xi = P(\xi)$.
Then $P$ has a coordinate representation $\td{P}(z,\zeta) = \paren{Az, N(z,\zeta)}$ in a neighborhood of $\xi$ where $A$ is an invertible matrix, $N(0,0) = 0$, and $\spec D_\zeta N(0,0) = \set{0}$.
Therefore given $\veps > 0$ we can apply Lemma~\ref{lem:superstable} to the nonlinearity $\td{P}(z,\zeta) - \paren{Az,0} = \paren{0,N(z,\zeta)}$ to find $\delta,C > 0$ such that for all $\paren{z,\zeta}\in B_\delta(0)$ and $k\in\N$:
\eqnal{
\norm{\td{P}^k(z,\zeta) - \paren{A^kz,0}} \le C\veps^k\norm{\paren{z,\zeta}}.
}
We conclude that $P$ is arbitrarily well--approximated near $\xi$ by a diffeomorphism on a submanifold whose dimension equals $\rk DP^n(\xi)$.
}

\iftoggle{sds}{
\sect{$C^1$ Linearization}\label{app:c1}
The technical result we desire was originally established by Hartman in the course of proving that an invertible contraction is $C^1$--conjugate to its linearization\footnote{Readers may be more familiar with the Hartman--Grobman Theorem \see{\cite[Theorem~1.4.1]{GuckenheimerHolmes1983} or \cite[Theorem~7.8]{Sastry1999}} which states that the phase portrait near an exponentially stable fixed point of a discrete--time dynamical system is \emph{topologically} conjugate to its linearization.}.
The original statement in~\cite{Hartman1960a} only considered invertible contractions.
However, as noted in~\cite{AulbachGaray1994}, the proof in~\cite{Hartman1960a} of the result we require does not make use of invertibility and the conclusion is still valid if zero is an eigenvalue of the linearization.
For details we refer the reader to \cite{Abbaci2004}, which also contains a generalization to \emph{hyperbolic} periodic orbits whose eigenvalues satisfy genericity conditions.


\begin{theorem}[Induction Assertion in~\cite{Hartman1960a}]\label{thm:hartman}
Let $U\subset\R^n$ be a neighborhood of the origin and $P:U\into\R^n$ a $C^1$ map of the form
\eqnal{
P(x,y,z) &= \paren{Ax + X(x,y,z),\ By + Y(x,y,z),\ Cz}
}
such that 
\eqnal{
DP(0) &= \mat{ccc}{A & 0 & 0 \\ 0 & B & 0 \\ 0 & 0 & C}
}
where: 
\begin{enumerate}
\item $x\in\R^k$, $y\in\R^\ell$, $z\in\R^m$ and $k+\ell+m = n$; 
\item $A\in\R^{k\times k}$, $B\in\R^{\ell\times \ell}$, and $C\in\R^{m\times m}$;
\item $X:\R^n\into\R^k$ and $Y:\R^n\into\R^\ell$ are $C^1$; 
\item $D_x X$, $D_y X$, $D_x Y$, and $D_y Y$ are uniformly Lipschitz continuous in $(x,y)$; 
\item $D_z X$ and $D_z Y$ are uniformly H\"{o}lder continuous in $z$;
\end{enumerate}
Suppose all the eigenvalues of $B$ have the same magnitude, that the eigenvalues of $A$ have smaller magnitude and those of $C$ have larger magnitude than those of $B$, and all eigenvalues of $DP(0)$ lie inside the unit disc:
\eqnal{\label{eqn:hartman}
\forall b,\beta\in\spec B &: \abs{b} = \abs{\beta};\\ 
\forall a\in\spec A, b\in\spec B, c\in\spec C &: 0 \le \abs{a} < \abs{b} < \abs{c} < 1.
}
Then there is a neighborhood of the origin $V\subset\R^n$ and a $C^1$ diffeomorphism $\vphi:V\into\R^n$ of the form
\eqnal{
\vphi(x,y,z) &= \paren{x + \vphi_X(z),\ y + \vphi_Y(x,y,z),\ z}
}
for which $D\vphi(0) = I$ and for all $(u,v,w)\in\vphi(V)$ we have 
\eqnal{
(\vphi\circ P\circ\vphi^{-1})(u,v,w) = \paren{Au + U(u,v,w),\ Bv,\ Cw}
}
where: 
\begin{enumerate}
\item $U:\vphi(V)\into\R^k$ is $C^1$;
\item $D_u U$ is uniformly Lipschitz continuous in $(u,v,w)$;
\item $D_v U$ and $D_w U$ are uniformly Lipschitz continuous in $u$;
\item $D_v U$ and $D_w U$ are uniformly H\"{o}lder continuous in $(v,w)$.
\end{enumerate}
\end{theorem}

\rem{Theorem~\ref{thm:hartman} may be applied inductively to exactly linearize all eigendirections corresponding to non--zero eigenvalues via a $C^1$ change--of--coordinates; this is the content of Lemma~\ref{lem:approx} in Section~\ref{sec:sds:approx}.}
}

\bibliographystyle{IEEEtran}
\bibliography{floq}

\end{document}